\documentclass{article}

\usepackage{geometry}
\usepackage[utf8]{inputenc} 
\usepackage{amsmath,amssymb}
\usepackage{tikz}
\usepackage{pgfplots} 
\usepackage{graphicx,subcaption}
\graphicspath{{./figs/}{./}{./Images/}}
\usetikzlibrary{patterns}
\usepackage{dsfont}
\usepackage{xcolor}
\usepackage{csquotes} 
\usepackage[backend=biber,style=numeric,sorting=none,language=french,maxbibnames=99]{biblatex}
\usepackage{url}            
\usepackage[T1]{fontenc}    
\usepackage[utf8]{inputenc} 

\usepackage{amsmath}        
\usepackage{mathrsfs}       
\usepackage{amssymb}        
\usepackage{amsfonts}       
\usepackage{cancel}

\usepackage{graphicx} 
\usepackage{wrapfig}  

\usepackage{tikz}     
\usepackage[framemethod=TikZ]{mdframed}

\usepackage{color}
\usepackage{dsfont}
\usepackage{bbm}
\usepackage{bm}

\usepackage{transparent} 

\usepackage[retain-zero-exponent=true]{siunitx} 

\PassOptionsToPackage{svgnames}{color}

\newcommand{\pa}{{\partial}\,} 

\newtheorem{rmq}{Remark}

\newtheorem{prop}{Proposition}
\newtheorem{prty}{Property}
\newtheorem{hypo}{Hypothesis}
\newcommand{\ds}{\displaystyle}

\newcommand{\ov}{\overline}

\DeclareMathOperator{\Deltavec}{\mathbf{\Delta}}

\newcommand{\dive}{\mathrm{div}}

\newcommand{\grad}{{\bf \ds \nabla}\,}

\newcommand{\divh}{{\rm div_h}\,}

\newcommand{\Th}{{\mathcal{T}_h}\,}

\def\eps{{\varepsilon}}

\def\Ical{\mathcal I}
\def\Fcal{\mathcal F}

\def\Mcal{\mathcal M}

\def\Pcal{\mathcal P}

\def\Scal{\mathcal S}
\def\Tcal{\mathcal T}

\def\bbL{{\mathbb L}}

\def\N{{\mathbb N}}

\def\R{{\mathbb R}}

\def\bH{{\bf H}}

\def\bL{{\bf L}}
\def\bP{{\bf P}}

\def\bV{{\bf V}}
\def\bX{{\bf X}}

\def\evec{{\bm e}}
\def\fvec{{\bm f}}

\def\nvec{{\bm n}}

\def\uvec{{\bm u}}
\def\vvec{{\bm v}}
\def\wvec{{\bm w}}
\def\xvec{{\bm x}}

\renewcommand{\tilde}{\widetilde}

\newcommand{\ul}{\underline}

\def\R{\mathbb{R}}
\def\bH{\mathbf{H}}
\def\bL{\mathbf{L}}
\def\bP{\mathbf{P}}

\def\bX{\mathbf{X}}
\def\fvec{\mathbf{f}}
\def\nvec{\mathbf{n}}

\def\uvec{\mathbf{u}}
\def\vvec{\mathbf{v}}
\def\wvec{\mathbf{w}}
\def\xvec{\mathbf{x}}

\def\Fcal{\mathcal{F}}
\def\Gcal{\mathcal{G}}
\def\Ical{\mathcal{I}}
\def\Mcal{\mathcal{M}}
\def\Scal{\mathcal{S}}
\def\Tcal{\mathcal{T}}
\def\Om{\Omega}
\def\Grad{\mathrm{\mathbf{Grad}}\,}
\def\grad{\mathrm{\mathbf{grad}}\,}

\def\dive{\mathrm{div}\,}

\newcommand{\logLogSlopeTriangle}[5]
{
	
	\pgfplotsextra
	{
		\pgfkeysgetvalue{/pgfplots/xmin}{\xmin}
		\pgfkeysgetvalue{/pgfplots/xmax}{\xmax}
		\pgfkeysgetvalue{/pgfplots/ymin}{\ymin}
		\pgfkeysgetvalue{/pgfplots/ymax}{\ymax}
		
		\pgfmathsetmacro{\xArel}{#1}
		\pgfmathsetmacro{\yArel}{#3}
		\pgfmathsetmacro{\xBrel}{#1-#2}
		\pgfmathsetmacro{\yBrel}{\yArel}
		\pgfmathsetmacro{\xCrel}{\xArel}
		
		\pgfmathsetmacro{\lnxB}{\xmin*(1-(#1-#2))+\xmax*(#1-#2)} 
		\pgfmathsetmacro{\lnxA}{\xmin*(1-#1)+\xmax*#1} 
		\pgfmathsetmacro{\lnyA}{\ymin*(1-#3)+\ymax*#3} 
		\pgfmathsetmacro{\lnyC}{\lnyA+#4*(\lnxA-\lnxB)}
		\pgfmathsetmacro{\yCrel}{\lnyC-\ymin)/(\ymax-\ymin)} 
		
		\coordinate (A) at (rel axis cs:\xArel,\yArel);
		\coordinate (B) at (rel axis cs:\xBrel,\yBrel);
		\coordinate (C) at (rel axis cs:\xCrel,\yCrel);
		
		\draw[#5]   (A)-- node[pos=0.5,anchor=north] {1}
		(B)-- 
		(C)-- node[pos=0.5,anchor=west] {#4}
		cycle;
	}
}

\definecolor{myblue}{rgb}{0.0, 0.0, 1.0}
\definecolor{myred}{rgb}{0.75, 0.0, 0.2}
\definecolor{mygreen}{rgb}{0.01, 0.75, 0.24}

\begin{document}

\title{Improved Crouzeix-Raviart scheme for the Stokes and Navier-Stokes problem}
\author{E. Chénier, \footnote{ Univ Paris Est Creteil,  CNRS, UMR 8208, MSME, F-77454 Marne-la-Vallée, France \newline \texttt{eric.chenier@univ-eiffel.fr}}
\and 
E. Jamelot, \footnote{Universit\'e Paris-Saclay, CEA,  Service de Thermo-hydraulique et de M\'ecanique des Fluides, 91191 Gif-sur-Yvette, France \texttt{erell.jamelot@cea.fr}}
\and 
C. Le Potier, \footnote{Universit\'e Paris-Saclay, CEA,  Service de Thermo-hydraulique et de M\'ecanique des Fluides,  91191 Gif-sur-Yvette, France \texttt{christophe.le-potier@cea.fr}}
\and
A. Peitavy \footnote{Universit\'e Paris-Saclay, CEA, Service de Thermo-hydraulique et de M\'ecanique des Fluides, 91191 Gif-sur-Yvette, France \texttt{andrew.peitavy@cea.fr} }

}

\newpage

\maketitle 

\begin{abstract}

The resolution of the incompressible Navier-Stokes equations is tricky, and it is well known that one of the major issue is to approach the space: $$ \bH^1(\Omega) \cap \bH(\dive0;\Omega):=\left\{\vvec\in \bH^1(\Omega)\, :\, \dive\vvec = 0 \right\}. $$ 
The non-conforming Crouzeix-Raviart finite element are convenient since they induce local mass conservation. Moreover they are such that the stability constant of the Fortin operator is equal to $1$. This implies that they can easily handle anisotropic mesh \cite{Apel01,Apel01Bis}. However spurious velocities may appear and damage the approximation.

 We propose a scheme here that allows to reduce the spurious velocities. It is based on a new discretisation for the gradient of pressure based on the symmetric MPFA scheme (finite volume MultiPoint Flux Approximation)  \cite{agelasMasson08,LEPOTIER05bis,lepotierhdr}.

\end{abstract}

\section{Motivation}\label{sec:intro}
The TrioCFD code is a computational fluid dynamics (CFD) simulation software
developed at the CEA. It is open source, object-oriented and massively parallel. It is
dedicated to the numerical simulation of turbulent flows for scientific and industrial
applications, particularly in the nuclear field. Let $\Om$, the domain of study, be an open connected bounded domain of $\R^d$, $d=2,\,3$, with a polygonal $(d=2)$ or Lipschitz polyhedral $(d=3)$ boundary $\pa\Om$ with constant physical properties. Let $T>0$ be a simulation time. The TrioCFD code solves the incompressible Navier-Stokes equations which read: Find $(\uvec(\xvec,t),p(\xvec,t))$ such that $\forall(\xvec,t)\in\Om\times(0,T)$,
\begin{equation}\label{eq:NS}
\left\{
\begin{array}{rcl}
\ds\pa_t\uvec-\nu\Deltavec\uvec+(\uvec\cdot\grad)\uvec+\grad p&=&\fvec,\\
\dive\uvec&=&0,\\
u(\xvec,0)&=&u_0(\xvec).
\end{array}\right.
\end{equation}
We consider here Dirichlet boundary conditions for the velocity $\uvec$ and we impose a normalization condition for the pressure $p$:
\[
\uvec=0\mbox{ on }\pa\Om,\quad\int_\Om p=0.
\]
The vector field $\uvec$ represents the velocity of the fluid and the scalar field $p$ represents its pressure divided by the fluid density which is supposed to be constant. 
The first equation of \eqref{eq:NS} corresponds to the momentum balance equation and the second one corresponds to the mass conservation. The constant parameter $\nu>0$ is the kinematic viscosity of the fluid. The vector field $\fvec$ represents the body force divided by the fluid density. We first consider the steady Stokes problem which reads:
\begin{equation}\label{eq:Stokes}
\mbox{Find }(\uvec,p)\mbox{  such that }\forall\xvec\in\Om:\,\left\{\begin{array}{rcl}
-\nu\Delta\uvec+\grad p&=&\fvec,\\
\dive\uvec&=&0.
\end{array}\right.
\end{equation}
Before stating the variational formulation of Problem \eqref{eq:Stokes}, we provide some definition and reminders. Let us set $\bL^2(\Om)=(L^2(\Om))^d$, $\bH^1_0(\Om)=(H^1_0(\Om))^d$, $\bH^{-1}(\Omega)=(H^{-1}(\Om))^d$ its dual space and $L^2_{zmv}(\Om)=\{q\in L^2(\Om)\,|\,\int_{\Om}q=0\}$. We recall that $\bH(\dive;\,\Om)=\{\vvec\in\bL^2(\Om)\,|\,\dive\vvec\in\, L^2(\Om)\}$. Let us first recall Poincaré-Steklov inequality:
\begin{equation}\label{eq:Poincare}
\exists C_{PS}>0\,|\,\forall v\in H^1_0(\Om),\quad \|v\|_{L^2(\Om)}\leq C_{PS}\|\grad v\|_{\bL^2(\Om)}.
\end{equation}
Thanks to this result, in $H^1_0(\Om)$, the semi-norm is equivalent to the natural norm, so that the scalar product reads $(v,w)_{H^1_0(\Om)}=(\grad v,\grad w)_{\bL^2(\Om)}$ and the norm is $\|v\|_{H^1_0(\Om)}=\|\grad v\|_{\bL^2(\Om)}$. Let $\vvec,\,\wvec\in\bH^1_0(\Om)$, we denote by $(v_i)_{i=1}^d$ (resp. $(w_i)_{i=1}^d$) the components of $\vvec$ (resp. $\wvec$), and we set $\Grad\vvec=(\pa_j v_i)_{i,j=1}^d\in\bbL^2(\Om)$, where $\bbL^2(\Om)=[L^2(\Om)]^{d\times d}$. We have: 
\[
(\Grad\vvec,\Grad\wvec)_{\bbL^2(\Om)}=(\vvec,\wvec)_{\bH^1_0(\Om)}=\ds\sum_{i=1}^d(v_i,w_i)_{H^1_0(\Om)}\mbox{ and }\|\vvec\|_{\bH^1_0(\Om)}=\|\Grad\vvec\|_{\bbL^2(\Om)}.\]
Let us set $\bV=\left\{\vvec\in\bH^1_0(\Om)\,|\,\dive\vvec=0\right\}$. The space $\bV$ is a closed subset of $\bH^1_0(\Om)$. We denote by $\bV^\perp$ the orthogonal of $\bV$ in $\bH^1_0(\Om)$. We recall that \cite[cor. I.2.4]{GiRa86}: 
\begin{prop}\label{prop:diviso}
	The operator $\dive:\,\bH^1_0(\Om)\rightarrow L^2(\Om)$ is an isomorphism of $\,\bV^{\perp}$ onto $L^2_{zmv}(\Om)$. We call $C_{\dive}$ the constant such that:
	\begin{equation}\label{eq:BAI}
	\forall p\in L^2_{zmv}(\Om),\,\exists!\vvec\in\bV^{\perp}\,|\,\dive\vvec=p\mbox{ and }
	\|\vvec\|_{\bH^1_0(\Om)}\leq C_{\dive}\|p\|_{L^2(\Om)}.
	\end{equation}
\end{prop}
Let us set~:
\begin{equation}
\label{eq:BilinForms}
a_\nu:\left\{\begin{array}{rcl}\bH^1_0(\Omega)\times\bH^1_0(\Omega)&\rightarrow&\R\\
(\uvec', \vvec)&\mapsto&\nu\,(\uvec', \vvec)_{\bH^1_0(\Om)}
\end{array}\right.\mbox{ and }b:\left\{\begin{array}{rcl}\bH^1_0(\Omega)\times L^2_{zmv}(\Om)&\rightarrow&\R\\
(\vvec, q)&\mapsto&(\dive\vvec, q)_{L^2(\Om)}
\end{array}\right..
\end{equation}
Classically, the variational formulation of Problem \eqref{eq:Stokes} reads: 
\begin{equation}
\label{eq:StokesVF}
\mbox{Find }(\uvec, p) \in \bH^1_0(\Omega)\times L_{zmv}^2(\Omega)\,|\,\left\{
\begin{array}{rcll}
a_\nu(\uvec, \vvec)_{\bH^1_0(\Om)}-b(\vvec,p)&=&\langle\fvec, \vvec\rangle & \forall \vvec \in  \bH^1_0(\Omega),\\
b(\uvec,q)&=&0 &\forall q \in L_{zmv}^2(\Omega).
\end{array}
\right.
\end{equation}
This saddle point problem is well-posed. Indeed, the bilinear form $a_\nu(\cdot,\cdot)$ is continuous and coercive. Moreover, the bilinear form $b(\cdot,\cdot)$ is continuous and due to Proposition \ref{prop:diviso}, it satisfies the following inf-sup condition:
\begin{equation}
\label{eq:CIS}
\forall q\in L^2_{zmv}(\Om)\backslash\{0\},\quad\exists\,\vvec_q\in\bH^1_0(\Omega)\backslash\{0\}\,|\quad\frac{b(\vvec_q,q)}{\|\vvec_q\|_{\bH^1_0(\Om)}\,\|q\|_{L^2(\Om)}}\geq C_{\dive}.
\end{equation}
In TrioCFD code, the spatial discretization of Problem \eqref{eq:Stokes} is based on  first order nonconforming Crouzeix-Raviart finite element method.\\
The outline of this article is as follows: in section \ref{sec:disc}, we provide some notations for the discretization. Next, in section \ref{sec:P1nc-P0}, we recall the first order nonconforming finite element method, that we call the $\bP^1_{nc}-P^0$ scheme. Then in section \ref{sec:P1nc-P0P1}, we describe the spatial discretization of TrioCFD code for simplicial meshes. We call this discretization the $\bP^1_{nc}-(P^0+P^1)$ scheme. This discretization is very precise in $2D$. It is also precise in $3D$, except when the source term is a strong gradient. In order to obtain the same accuracy in $3D$ than in $2D$, one must increase the number of degrees of freedom of the discrete pressure space, which leads to a more expensive numerical scheme. Our aim is to develop a new numerical scheme that would be precise both in $2D$ and $3D$, but at a lower cost. We present such a scheme in section \ref{sec:MPFA} and numerical illustration in section \ref{sec:numstokes} and \ref{sec:numNstokes}.
\section{Discrete notations}\label{sec:disc}
We call $(O,(x_{d'})_{d'=1}^d)$ the Cartesian coordinates system, of orthonormal basis $(\evec_{d'})_{d'=1}^d$. 
Consider $(\Tcal_h)_h$ a simplicial triangulation sequence of $\Om$. For a triangulation $\Tcal_h$, we use the following index sets:
\begin{itemize}
	\item $\Ical_K$ denotes the index set of the elements, such that $\Tcal_h:=\ds\bigcup_{\ell\in\Ical_K}K_\ell$  is the set of elements.
	\item $\Ical_F$ denotes the index set of the facets\footnote{The term facet stands for face (resp. edge) when $d=3$ (resp. $d=2$).}, such that $\Fcal_h:=\ds\bigcup_{f\in\Ical_F}F_f$ is the set of facets.
	\item[]Let $\Ical_F=\Ical_F^i\cup\Ical_F^b$, where $\forall f\in\Ical_F^i$, $F_f\in\Om$ and $\forall f\in\Ical_F^b$, $F_f\in\pa\Om$.
	\item $\Ical_S$ denotes the index set of the vertices, such that $(S_j)_{j\in\Ical_S}$ is the set of vertices.
	\item[]Let $\Ical_S=\Ical_S^i\cup\Ical_S^b$, where $\forall j\in\Ical_S^i$, $S_j\in\Om$ and $\forall j\in\Ical_S^b$, $S_j\in\pa\Om$.
\end{itemize} 
We also define the following index subsets:
\begin{itemize}
	\item $\forall\ell\in\Ical_K$, $\Ical_{F,\ell}=\{f\in\Ical_F\,|\,F_f\in K_\ell\},\quad\Ical_{S,\ell}=\{j\in\Ical_S\,|\,S_j\in K_\ell\}$.
	\item $\forall j\in\Ical_S$, $\Ical_{K,j}=\{\ell\in\Ical_K\,|\,S_j\in K_\ell\}$, $\quad N_{K,j}:=\mathrm{card}(\Ical_{K,j})$.
	\item $\forall j\in\Ical_S$, $\Ical_{S,j}=\{k\in\Ical_S\,|\,S_kS_i\in\Fcal_h\}$, $\quad N_{S,j}:=\mathrm{card}(\Ical_{S,j})$.
\end{itemize}
Notice that in $2D$, $N_{K,j}=N_{S,j}$.\\
For all $f\in\Ical_F$, $M_f$ denotes the barycentre of $F_f$, and by $\nvec_f$ a unit normal (outward oriented if $F_f\in\pa\Omega$). For all $j\in\Ical_S$, for all $\ell\in\Ical_{K,j}$, $\lambda_{j,\ell}$ denotes the barycentric coordinate of $S_j$ in $K_\ell$; $F_{j,\ell}$ denotes the face opposite to vertex $S_j$ in element $K_\ell$. We call $\Scal_{j,\ell}$ the outward normal vector of $F_{j,\ell}$ and of norm $|\Scal_{j,\ell}|=|F_{j,\ell}|$. Let introduce spaces of piecewise regular elements:\\
We set $\Pcal_h H^1=\left\{v\in L^2(\Om)\,;\quad\forall \ell\in\Ical_K,\,v_{|K_\ell}\in H^1(K_\ell)\right\}$, endowed with the scalar product~:
\[
(v,w)_h:=\sum_{\ell\in\Ical_K}(\grad v,\grad w)_{\bL^2(K_\ell)}\quad
\|v\|_h^2=\sum_{\ell\in\Ical_K}\|\grad v\|^2_{\bL^2(K_\ell)}.
\]
We set $\Pcal_h\bH^1=[\Pcal_hH^1]^d$, endowed with the scalar product~:
\[
(\vvec,\wvec)_h:=\sum_{\ell\in\Ical_K}(\Grad\vvec,\Grad\wvec)_{\bbL^2(K_\ell)}\quad
\|\vvec\|_h^2=\sum_{\ell\in\Ical_K}\|\Grad \vvec\|^2_{\bbL^2(K_\ell)}.
\]
Let $f\in\Ical_F^i$ such that $F_f=\pa K_L\cap\pa K_R$ and let $\nvec_f$ the unit normal that is outward $K_L$ oriented.\\
The jump (resp. average) of a function $v\in \Pcal_h H^1$ across the facet $F_f$, in $\nvec_f$ direction, is defined as follows: $[v]_{F_f}:=v_{|K_L}-v_{|K_R}$ (resp. $\{v\}_{F_f}:=\frac{1}{2}(v_{|K_L}+v_{|K_R})\,$). For $f\in\Ical_F^b$, we set: $[v]_{F_f}:=v_{|F_f}$ and $\{v\}_{F_f}:=v_{|F_f}$.\\
We set $\Pcal_h\bH(\dive)=\left\{\vvec\in \bL^2(\Omega)\,;\quad\forall\ell\in\Ical_K,\,\vvec_{|K_\ell}\in\bH(\dive;\,K_\ell)\right\}$, and we define the operator $\divh$ such that: 
\[
\forall\vvec\in\Pcal_h\bH(\dive),\,\forall q\in L^2(\Om),\quad
(\divh\vvec,q)=\sum_{\ell\in\Ical_K}(\dive\vvec,q)_{L^2(K_\ell)}.
\]
For all $D\subset\R^d$, and $k\in\N^*$, we call $P^k(D)$ the set of order $k$ polynomials on $D$, $\bP^k(D)=(P^k(D))^d$, and we consider the space of the broken polynomials:
\[
P^k_{disc}(\Tcal_h)=\left\{q\in L^2(\Om);\quad \forall\ell\in\Ical_K,\,q_{|K_\ell}\in P^k(K_\ell)\right\},\quad\bP^k_{disc}(\Tcal_h):=(P^k_{disc}(\Tcal_h))^d.
\]
We let $P^0(\Tcal_h)$ be the space of piecewise constant functions on $\Tcal_h$.
\begin{equation}
\label{eq:Qhk}
\forall k \in \mathbb{N}, \quad Q_{k,h}:= P^k(\Th) \cap L^2_{zmv}(\Om)
\end{equation}
We will now describe three numerical scheme to solve \eqref{eq:Stokes} for which the components of the velocity is discretized with the first order nonconforming Crouzeix-Raviart finite element method \cite[\S 5, Example 4]{CrRa73}. For simplicity, we suppose now that $\fvec\in\bL^2(\Om)$.
\section{The $\bP^1_{nc}-P^0$ scheme}\label{sec:P1nc-P0}
The first order nonconforming finite element method was introduced by Crouzeix and Raviart in the seminal paper \cite{CrRa73} to solve Stokes Problem \eqref{eq:Stokes}. We call it the $\bP^1_{nc}-P^0$ scheme. Let us consider $X_{h}$ (resp. $X_{0,h}$), the space of nonconforming approximation of $H^1(\Om)$ (resp. $H^1_0(\Om)$) of order $1$:
\begin{equation}
\label{eq:Xh}
X_{h}=\left\{v_h\in P^1_{disc}(\Tcal_h)\,;\quad\forall f\in\Ical_F^i,\,\int_{F_f} [v_h]=0\right\}.
\end{equation}
\begin{equation}
\label{eq:X0h}
X_{0,h}=\left\{v_h\in X_{h} \,;\quad\forall f\in\Ical_F^b,\,\int_{F_f} [v_h]=0\right\}.
\end{equation}
\begin{prop}\label{pro:broknorm}
	The broken norm $v_h\rightarrow\|v_h\|_h$ is a norm over $X_{0,h}$. 
\end{prop}
The following discrete Poincaré–Steklov inequality holds \cite[Lemma 36.6]{ErGu21-II}: there exists a constant $C_{PS}^{nc}>0$ such that
\begin{equation}\label{eq:cps-constant}
\forall v_h\in X_{0,h},\quad\|v_h\|_{L^2(\Om)}\leq C_{PS}^{nc}\,\|v_h\|_h.
\end{equation}
The constant $C_{PS}^{nc}$ is independent of the triangulation $\Tcal_h$ and it is proportional to the diameter of $\Om$.
\\
We can endow $X_{0,h}$ with the basis $(\psi_f)_{f\in\Ical_F^i}$ such that: $\forall \ell\in\Ical_K$,
\[
\psi_{f|K_\ell}=\left\{\begin{array}{cl}1-d\lambda_{i,\ell}&\mbox{if }f\in\Ical_{F,\ell},\\ 0&\mbox{ otherwise,}\end{array}\right.
\]
where $S_i$ is the vertex opposite to $F_f$ in $K_\ell$.  We then have $\psi_{f|F_f}=1$, so that $[\psi_f]_{F_f}=0$ if $f\in\Ical_F^i$ (i.e. $F_f\in \Om$), and $\forall f'\neq f$, $\int_{F_{f'}}\psi_f=0$. We have: $X_{0,h}=\mathrm{vect}\left((\psi_f)_{f\in\Ical_F^i}\right)$.\\ 
The Crouzeix-Raviart interpolation operator $\pi_h$ for scalar functions is defined by:
\[
\pi_h: \left\{\begin{array}{rcl} H^1(\Om)&\rightarrow&X_h\\ v&\mapsto&\ds\sum_{f\in\Ical_F}\pi_fv\,\psi_{f}\end{array}\right.,\mbox { where }\pi_fv=\frac{1}{|F_f|}\int_{F_f} v.
\]
Notice that $\forall f\in\Ical_F$, $\int_{F_f} \pi_hv=\int_{F_f} v$. Moreover, the Crouzeix-Raviart interpolation operator preserves the constants, so that $\pi_h\ul{v}_\Om=\ul{v}_\Om$ where $\ul{v}_\Om=\int_\Om v/|\Om|$. 

The space of nonconforming approximation $\bH^1_0(\Om)$ of order $1$ is $\bX_{0,h}=(X_{0,h})^d$. For a vector $\vvec\in\bH^1(\Om)$ of components $(v_{d'})_{d'=1}^d$, the Crouzeix-Raviart interpolation operator is such that: $\Pi_h\vvec=\left(\pi_hv_{d'}\right)_{d'=1}^d$. We recall the following result:
\begin{prop}\label{lem:PiCR}
	The Crouzeix-Raviart interpolation operator $\Pi_h$ can play the role of the Fortin operator:
	\begin{align}
	\label{eq:WellPosed-CR-1}
	\forall\vvec\in\bH^1(\Om)\hspace*{5mm}\|\Pi_h\vvec\|_h&\leq\|\Grad\vvec\|_{\bbL^2(\Om)},\\
	\label{eq:WellPosed-CR-2}
	\forall\vvec\in\bH^1(\Om)\quad(\divh\Pi_h\vvec,q_h)&=(\dive\vvec,q_h)_{L^2(\Om)},\quad\forall q\in Q_h.
	\end{align}
	Moreover, for all $\vvec\in\bP^1(\Om)$, $\Pi_h\vvec=\vvec$.
\end{prop}
Notice that the stability constant of the bound on $\|\Pi_h\vvec\|_h$ is equal to $1$ \cite[Lemma 2]{Apel01}~: it is independent of the mesh.\\
Let us set $Q_h= Q_{0,h}$. We now define the following bilinear forms~:
\begin{equation}
\label{eq:DiscBilinForms}
a_{\nu,h}:\left\{\begin{array}{rcl}\bX_{0,h}\times\bX_{0,h}&\rightarrow&\R\\
(\uvec'_h, \vvec_h)&\mapsto&\nu\,(\uvec'_h, \vvec_h)_h
\end{array}\right.\mbox{ and }b_h:\left\{\begin{array}{rcl}\bX_{0,h}\times Q_h&\rightarrow&\R\\
(\vvec_h, q_h)&\mapsto&-(\divh\vvec_h, q_h)
\end{array}\right..
\end{equation}
We suppose here that $\fvec\in\bL^2(\Om)$. The discretization of variational formulation \eqref{eq:StokesVF} reads:
\begin{equation}
\label{eq:StokesVFh}
\mbox{Find }(\uvec_h, p_h) \in\bX_{0,h}\times Q_h\,|\,\left\{
\begin{array}{rcll}
a_{\nu,h}(\uvec_h, \vvec_h)_h+b_h(\vvec_h,p)&=&(\fvec, \vvec_h)_{L^2(\Om)} & \forall \vvec_h \in  \bX_{0,h},\\
b_h(\uvec_h,q_h)&=&0 &\forall q_h \in Q_h.
\end{array}
\right.
\end{equation}
This saddle point problem is well-posed. Indeed, the bilinear form $a_{\nu,h}(\cdot,\cdot)$ is continuous and coercive. Moreover, the bilinear form $b_h(\cdot,\cdot)$ is continuous and due to Proposition \ref{prop:diviso} and Lemma \ref{lem:PiCR}, it satisfies the following discrete inf-sup condition:
\begin{equation}
\label{eq:CISh}
\forall q_h\in Q_h\backslash\{0\},\quad\exists\,\vvec_h,q\in\bX_{0,h}\backslash\{0\}\,|\quad\frac{b_h(\vvec_h,q_h)}{\|\vvec_h\|_h\,\|q_h\|_{L^2(\Om)}}\geq C_{\dive}.
\end{equation}
Suppose that it exists $\phi\in H^1(\Om)\cap L^2_{zmv}(\Om)$ such that $\fvec=\grad\phi$. In that case, the solution to Problem \eqref{eq:Stokes} is $(\uvec,p)=(0,\phi)$. By integrating by parts, we have:
\[
\forall\vvec_h\in\bX_{0,h},\quad
(\fvec, \vvec_h)_{L^2(\Om)}=-(\divh\vvec_h,\phi)+\sum_{f\in\Ical_F^i}\int_{F_f}[\vvec_h\cdot\nvec_f]\,\phi.
\]
The term with the jump acts as a numerical source, which numerical influence is proportional to $1/\nu$. Hence, we cannot obtain exactly $\uvec_h=0$. There are different strategies to cure this well-known problem:
\begin{itemize}
	\item Using a polynomial approximation of higher degree \cite{FoSo83}.
	\item Increasing the space of the discrete pressures \cite{Heib03,Fort06}.
	\item Projecting the test-function on a discrete subspace of $\bH(\dive;\Om)$ \cite{Link14}.
\end{itemize}
We propose in the next section to give details on the second strategy.
\section{The $\bP^1_{nc}-(P^0+P^1)$ scheme}\label{sec:P1nc-P0P1}
In his thesis \cite{Heib03}, Heib proposed to use the following space discrete pressures space (cf. \eqref{eq:Qhk}):
\begin{equation}\label{eq:Qh}
\tilde{Q}_h=Q_{0,h}\oplus Q_{1,h}.
\end{equation}
For any $\tilde{q}_h\in\tilde{Q}_h$, we write: $\tilde{q}_h=q_{0,h}+q_{1,h}$, where $q_{0,h}\in Q_{0,h}$ and $q_{1,h}\in Q_{1,h}$. \\
Let us consider the following bilinear form:
\begin{equation}
\label{eq:btilde}
\tilde{b}_h:\left\{
\begin{array}{rcl}\bX\times \tilde{Q}_h&\rightarrow&\R \\ 
(\vvec_h,\tilde{q}_h)&\mapsto&-(\divh\vvec_h,q_{0,h})+(\vvec_h,\grad q_{1,h})_{\bL^2(\Om)}
\end{array}
\right..
\end{equation}
The discretization of Problem \eqref{eq:StokesVF} with $\bP^1_{nc}-(P^0+P^1)$ finite elements reads:
\begin{equation}
\label{eq:StokesVFhTrio}
\mbox{Find }(\uvec_h, p_h) \in\bX_{0,h}\times \tilde{Q}_h\,|\,\left\{
\begin{array}{rcll}
a_{\nu,h}(\uvec_h, \vvec_h)_h+\tilde{b}_h(\vvec_h,p_h)&=&(\fvec, \vvec_h)_{L^2(\Om)} & \forall \vvec_h \in  \bX_{0,h},\\
\tilde{b}_h(\uvec_h,\tilde{q}_h)&=&0 &\forall \tilde{q}_h \in \tilde{Q}_h.
\end{array}
\right.
\end{equation}
We will need the following Hypothesis \cite[Hyp.~4.1]{BeHe00}:
\begin{hypo}\label{hyp:cis}
	We suppose that the triangulation $\Tcal_h$ is such that the boundary $\pa\Om$ contains at most one edge in dimension $d=2$ and at most two faces in dimension $d=3$, of the same element $K_\ell$, $\ell\in\Ical_K$.
\end{hypo}
Under Hypothesis \eqref{hyp:cis}, one can prove that the bilinear form $\tilde{b}_h(\cdot,\cdot)$ is continuous that it satisfies the following discrete inf-sup condition \cite[\S 4.2]{Heib03}:
\begin{equation}
\label{eq:CIShTrio}
\forall \tilde{q}_h\in \tilde{Q}_h\backslash\{0\},\quad\exists\,\vvec_h\in\bX_{0,h}\backslash\{0\}\,|\quad\frac{\tilde{b}_h(\vvec_{h,\tilde{q}_h},\tilde{q}_h)}{\|\vvec_{h,\tilde{q}_h}\|_h\,\|\tilde{q}_h\|_{L^2(\Om)}}\geq \tilde{C}_{\dive},
\end{equation}
where the constant $\tilde{C}_{\dive}$ is independent of the mesh size. Compared to $\bP^1_{nc}-P^0$ scheme, the $\bP^1_{nc}-(P^1+P^0)$ scheme gives a better approximation of the velocity in the sense that the discrete mass conservation equation is strengthened. Indeed, one can show, for $d=2$ that \cite[Theorem 4.3.2]{Fort06}:
\begin{prty}\label{thm:P1nc-P0P1}
	Let $\vvec_h\in\bV_h:=\{\wvec_h\in\bX_h\,|\quad\forall q_h\in\tilde{Q}_h,\quad\tilde{b}_h(\wvec_h,q_h)=0\}$. \\
	Then for $d=2$, we have: for all $q_{2,h}\in Q_{2,h}$, $\tilde{b}_h(\wvec_h,q_{2,h})=(\grad q_{2,h},\vvec_h)_{\bL^2(\Om)}=0$.
\end{prty}
The proof of Property \ref{thm:P1nc-P0P1} relies on a quadrature formula which uses the degrees of freedom of the discrete pressure. As this formula cannot be extended in $3D$, this property does not hold. To recover Property \ref{thm:P1nc-P0P1} in $3D$, we must introduce $P^2$ discrete pressure degrees of freedom, located on the edges of the mesh, as detailed in \cite{Fort06}. This increases the number of unknowns by the number of cells, which leads to an expensive linear system. Hence, we look for a numerical scheme which could be as precise in $3D$ than in $2D$, but at a lower cost. In the next section, we propose a new strategy, which relies on the multi-points flux approximation to discretize the pressure gradient term in \eqref{eq:Stokes}.

\section{The $\bP^1_{nc}-P^0_{Mps}$ scheme}\label{sec:MPFA}

Here, we use the symmetric MPFA scheme (where MPFA stands for {\em multi-points flux approximation}) to discretize the pressure gradient term in \eqref{eq:Stokes}, in the case of a simplicial mesh. The discrete pressure space remains $Q_h=Q_{0,h}$. We call this new scheme the $\bP^1_{nc}-P^0_{Mps}$ scheme.
\par
Let us consider the $2D$ case.
To design the scheme, as it has been initially done in \cite{LEPOTIER05bis, LEPOTIER05}, we start by splitting the triangles into three quadrangles, connecting the barycentre of the triangle to the midpoint of each edges. Considering some $ q_h \in Q_h$, we will calculate an affine approximation of $Q_h$ on each quadrangles. To do so,we need to add temporary auxiliary unknowns located on the third of the edges.
\\
Let us introduce some notations. \\

Let  $j \in \Ical_{S}$. We define the macro-element $\Mcal_j$ such that $\ds\overline{\Mcal}_j:=\bigcap_{\ell\in\Ical_{K,j}}\overline{K}_\ell$. Let renumber the vertices so that: $S_0=S_j$, $\Ical_{S,0}=\{1,\cdots,N_{S,0}\}$ and for all $i\in\Ical_{S,0}$, $S_iS_{i+1}\in\Fcal_h$ (setting $S_{N_{S,0}+1}=S_{N_{S,0}}$). For $i\in\Ical_{S,0}$ we denote by:
\begin{itemize}
	\item $K_i$ the triangle of vertices $S_0 S_i S_{i+1}$, and we call its barycentre $G_i$.
	\item $F_i$ the edge such that $F_i=S_0S_i$, and we call $M_i$ its midpoint. 
	\item $F_{i,0}$ the edge opposite to $S_0$ in $K_i$. 
	\item $\tilde{F}_i$ the half-edges defined by $S_0$ and the midpoint of $F_i$.
	\item $Q_i$ the quadrangle of vertices $S_0\,M_i\,G_i\,M_{i+1}$ (Fig. \ref{fig:Qi}-(a) for $S_0\in\Om$ and \ref{fig:Qj}-(a) for $S_0\in\pa\Om$).
\end{itemize}
For $i,\,j\in\Ical_{S,0}$, we denote by $\Scal_{i,j}$ the normal vector outgoing of $K_j$ at $F_i$ and of norm $|F_i|$. For $i\in\Ical_{S,0}$, we call $\Scal_{0,i}$ the normal vector outgoing of $K_i$ at $F_{i,0}$. On Figure \ref{fig:MS_and_T1}-(a), we represent $\Mcal_0$ in case $S_0\in\Om$ and $N_{S,0}=6$. On Figure \ref{fig:MS_and_T1}-(b), we represent the triangle $K_1$ with the vectors $(\Scal_{j,1})_{j=0}^2$ and its barycentre $G_1$.
\vspace{-0.3cm}
\begin{center}
	\begin{figure}[ht!]
		\begin{subfigure}[b]{0.39\linewidth}
			\begin{tikzpicture}[scale=0.5]
			\draw[-] (0,0) -- (3,4) -- (-2,4) --cycle;
			\draw[-] (0,0) -- (-4,-1) -- (0,-4) --cycle;
			\draw[-] (0,0) -- (3,-2) -- (4,1) --cycle;
			\draw[-] (-4,-1) -- (-2,4);
			\draw[-] (0,-4) -- (3,-2);
			\draw[-] (4,1) -- (3,4);
			\draw[dashed,orange] (1.5,2)--(1/3,8/3)--(-1,2)--(-2,1)--(-2,-0.5)--(-4/3,-5/3)--(0,-2)--(1,-2)--(1.5,-1)--(7/3,-1/3)--(2,0.5)--(7/3,5/3)--(1.5,2);
			\draw[dashed,orange] (1/3,8/3)--(1/2,8/2);
			\draw[dashed,orange] (-2,1)--(-3,3/2);
			\draw[dashed,orange] (-4/3,-5/3)--(-4/2,-5/2);
			\draw[dashed,orange] (1,-2)--(3/2,-3);
			\draw[dashed,orange] (7/3,-1/3)--(7/2,-1/2);
			\draw[dashed,orange] (7/3,5/3)--(7/2,5/2);
			
			\coordinate [label=left:$S_0$] (S_0) at (0.7,1);
			\coordinate [label=right:$S_1$] (S1) at (3,4);
			\coordinate [label=left:$S_2$] (S2) at (-2,4);
			\coordinate [label=left:$S_3$] (S3) at (-4,-1);
			\coordinate [label=below:$S_4$] (S4) at (0,-4);
			\coordinate [label=right:$S_5$] (S5) at (3,-2);
			\coordinate [label=right:$S_6$] (S6) at (4,1);
			\coordinate [label=above:\textcolor{teal}{$F_{1,0}$}] (F1) at (1/2,4);
			\coordinate [label=left:\textcolor{teal}{$F_{2,0}$}] (F1) at (-3,1.5);
			\coordinate [label=left:\textcolor{teal}{$F_{3,0}$}] (F1) at (-2,-2.5);
			\coordinate [label=right:\textcolor{teal}{$F_{4,0}$}] (F1) at (1.6,-3);
			\coordinate [label=right:\textcolor{teal}{$F_{5,0}$}] (F1) at (3.5,-0.5);
			\coordinate [label=right:\textcolor{teal}{$F_{6,0}$}] (F1) at (3.5,2.5);
			\coordinate [label=left:\textcolor{teal}{$F_1$}] (F1) at (2,3);
			\coordinate [label = left:\textcolor{teal}{$F_2$}] (F2) at (-1.3,2);
			\coordinate [label=below:\textcolor{teal}{$F_3$}] (F3) at (-2.5,-0.8);
			\coordinate [label=right:\textcolor{teal}{$F_4$}] (F4) at (0,-2.6);
			\coordinate [label=right:\textcolor{teal}{$F_5$}] (F5) at (1.7,-1);
			\coordinate [label=above:\textcolor{teal}{$F_6$}] (F6) at (2.6,0.6);
			\coordinate [label=below:\textcolor{red}{$K_1$}] (T1) at (-1/3,11/3);
			\coordinate [label=right:\textcolor{red}{$K_2$}] (T2) at (-3.4,0.2);
			\coordinate [label=right:\textcolor{red}{$K_3$}] (T3) at (-5/3,-7/3);
			\coordinate [label=above:\textcolor{red}{$K_4$}] (T4) at (2,-2.5);
			\coordinate [label=right:\textcolor{red}{$K_5$}] (T5) at (7/3,-0);
			\coordinate [label=right:\textcolor{red}{$K_6$}] (T6) at (6/3,2.5);
			\end{tikzpicture}
			\caption{Macro-element $\Mcal_0=S_1\,S_2\,S_3\,S_4\,S_5\,S_6$.}
		\end{subfigure}
		\hspace{0.2\linewidth}
		\begin{subfigure}[b]{0.39\linewidth}
			\begin{tikzpicture}[scale=0.5]
			=
			\draw[-] (0,0) -- (3,4) -- (-2,4) --cycle;
			\coordinate [label=below:$S_0$] (S) at (0,0);
			\coordinate [label=right:$S_1$] (S1) at (3,4);
			\coordinate [label=left:$S_2$] (S2) at (-2,4);
			\draw[dashed,orange] (1/3,8/3)--(1.5,2);
			\draw[dashed,orange] (1/3,8/3)--(-1,2.0);
			\draw[dashed,orange] (1/3,8/3)--(1/2,4);
			\fill[fill=blue] (-1,2) circle (0.1cm);
			\coordinate [label=left:\textcolor{blue}{$M_2$}] (M2) at (-1,2.2);
			\fill[fill=blue] (1.5,2) circle (0.1cm);
			\coordinate [label=right:\textcolor{blue}{$M_1$}] (M1) at (1.7,2.2);
			\coordinate [label = above:$G_1$] (G) at (-1/5,8/3);
			\fill[fill=red] (1/3,8/3) circle (0.12cm);
			\draw [->] (1.5,2) --(5.5,-1);
			\coordinate [label = right:$\Scal_{1,1}$] (N1) at (3.8,0.5);
			\draw [->] (-1,2) -- (-5,0);
			\coordinate [label = left:$\Scal_{2,1}$] (N2) at (-3.5,1.2);
			\draw [->] (0.5,4) -- (0.5,8);
			\coordinate [label = right:$\Scal_{0,1}$] (N0) at (0.5,6);
			\end{tikzpicture}
			\caption{Triangle $K_1=S_0\,S_1\,S_2$.}
		\end{subfigure}
		\caption{Notations in case $N_{S,0}=6$ and $j\in \Ical_S^i$.}
		\label{fig:MS_and_T1}
	\end{figure}
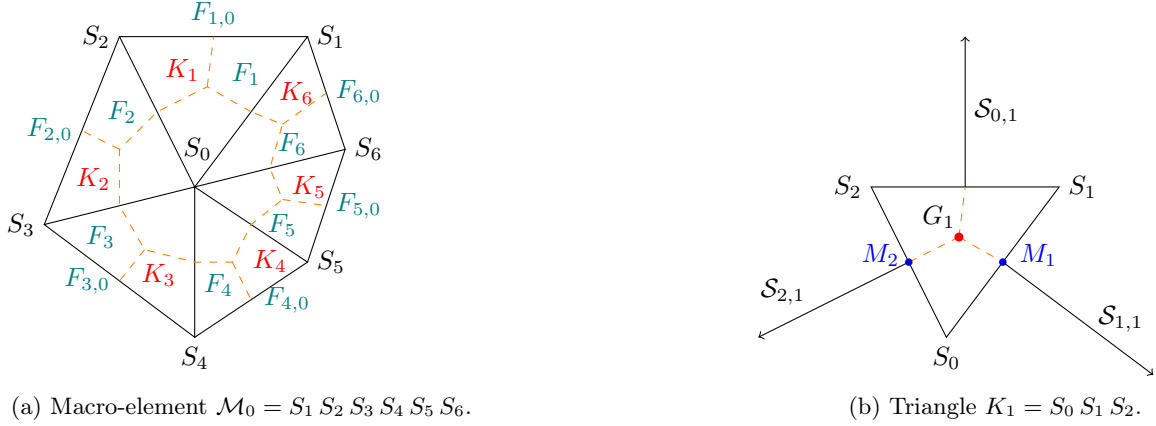
\end{center}
\vspace{-1.1cm}
\begin{center}
	\begin{figure}[ht!]
		\begin{subfigure}[b]{0.39\linewidth}
			\begin{tikzpicture}[scale=0.5]
			\draw[-] (0,0) -- (3,4) -- (-2,4) --cycle;
			\draw[-] (0,0) -- (-4,-1) -- (0,-4) --cycle;
			\draw[-] (0,0) -- (3,-2) -- (4,1) --cycle;
			\draw[-] (-4,-1) -- (-2,4);
			\draw[-] (0,-4) -- (3,-2);
			\draw[-] (4,1) -- (3,4);
			\draw[dashed,orange] (1.5,2)--(1/3,8/3)--(-1,2)--(-2,1)--(-2,-0.5)--(-4/3,-5/3)--(0,-2)--(1,-2)--(1.5,-1)--(7/3,-1/3)--(2,0.5)--(7/3,5/3)--(1.5,2);
			\coordinate [label=left:$S_0$] (S) at (0.7,1);
			\coordinate [label=right:$S_1$] (S1) at (3,4);
			\coordinate [label=left:$S_2$] (S2) at (-2,4);
			\coordinate [label=left:$S_3$] (S3) at (-4,-1);
			\coordinate [label=below:$S_4$] (S4) at (0,-4);
			\coordinate [label=right:$S_5$] (S5) at (3,-2);
			\coordinate [label=right:$S_6$] (S6) at (4,1);
			%
			\coordinate [label=below:\textcolor{orange}{$Q_1$}] (Q1) at (1/3,7.8/3);
			\coordinate [label=right:\textcolor{orange}{$Q_2$}] (Q2) at (-2,0.9);
			\coordinate [label=right:\textcolor{orange}{$Q_3$}] (Q3) at (-4/3,-4/3);
			\coordinate [label=above:\textcolor{orange}{$Q_4$}] (Q4) at (0.82,-1.9);
			\coordinate [label=right:\textcolor{orange}{$Q_5$}] (Q5) at (0.80,-0.30);
			\coordinate [label=right:\textcolor{orange}{$Q_6$}] (Q6) at (0.90,1);
			\end{tikzpicture}
			\caption{Quadrangles $(Q_i)_{i=1}^{N_{S,0}}$.}
		\end{subfigure}
		\hspace{0.2\linewidth}
		\begin{subfigure}[b]{0.39\linewidth} 
			\begin{tikzpicture}[scale=0.5]
			\draw[-] (0,0) -- (3,4) -- (-2,4) --cycle;
			\draw[-] (0,0) -- (-4,-1) -- (0,-4) --cycle;
			\draw[-] (0,0) -- (3,-2) -- (4,1) --cycle;
			\draw[-] (-4,-1) -- (-2,4);
			\draw[-] (0,-4) -- (3,-2);
			\draw[-] (4,1) -- (3,4);
			\draw[dashed,orange] (1.5,2)--(1/3,8/3)--(-1,2)--(-2,1)--(-2,-0.5)--(-4/3,-5/3)--(0,-2)--(1,-2)--(1.5,-1)--(7/3,-1/3)--(2,0.5)--(7/3,5/3)--(1.5,2);
			\coordinate [label=left:$S_0$] (S) at (0.7,1);
			\coordinate [label=right:$S_1$] (S1) at (3,4);
			\coordinate [label=left:$S_2$] (S2) at (-2,4);
			\coordinate [label=left:$S_3$] (S3) at (-4,-1);
			\coordinate [label=below:$S_4$] (S4) at (0,-4);
			\coordinate [label=right:$S_5$] (S5) at (3,-2);
			\coordinate [label=right:$S_6$] (S6) at (4,1);
			%
			\coordinate [label = above:$\overline{q}_1$] (pT1) at (1/3,8/3);
			\fill[fill=red] (1/3,8/3) circle (0.12cm);
			\coordinate [label = left:$\overline{q}_2$] (pT2) at (-2,1);
			\fill[fill=red] (-2,1) circle (0.12cm);
			\coordinate [label = left:$\overline{q}_3$] (pT3) at (-4/3,-5/3);
			\fill[fill=red] (-4/3,-5/3) circle (0.12cm);
			\coordinate [label = below:$\overline{q}_4$] (pT4) at (1,-2) ;
			\fill[fill=red] (1,-2) circle (0.12cm);
			\coordinate [label = right:$\overline{q}_5$] (pT5) at (7/3,-1/3) ;
			\fill[fill=red] (7/3,-1/3) circle (0.12cm);
			\coordinate [label = right:$\overline{q}_6$] (pT6) at (7/3,5/3) ;
			\fill[fill=red] (7/3,5/3) circle (0.12cm);
			\coordinate [label = above:$\tilde{q}_1$] (p1) at (0.9,4/3);
			\fill[fill=magenta] (1,4/3) circle (0.1cm);
			\coordinate [label = left:$\tilde{q}_2$] (p2) at (-2/3,1.2);
			\fill[fill=magenta] (-2/3,4/3) circle (0.1cm);
			\coordinate [label = below:$\tilde{q}_3$] (p3) at (-4/3,-1/3);
			\fill[fill=magenta] (-4/3,-1/3) circle (0.1cm);
			\coordinate [label = right:$\tilde{q}_4$] (p4) at (0,-4/3);
			\fill[fill=magenta] (0,-4/3) circle (0.1cm);
			\coordinate [label = above:$\tilde{q}_5$] (p5) at (1.3,-0.68);
			\fill[fill=magenta] (1,-2/3) circle (0.1cm);
			\coordinate [label = above:$\tilde{q}_6$] (p6) at (4/3,1/3);
			\fill[fill=magenta] (4/3,1/3) circle (0.1cm);
			
			\end{tikzpicture}
			\caption{Discrete pressures ($\overline{q}_i$,$\tilde{q}_i)_{i=1}^{N_{S,0}}$.}
		\end{subfigure}
		\caption{MPFA Scheme for $j \in \Ical_S^i$ and $N_{S,0}=6$.}
		\label{fig:Qi}
	\end{figure}
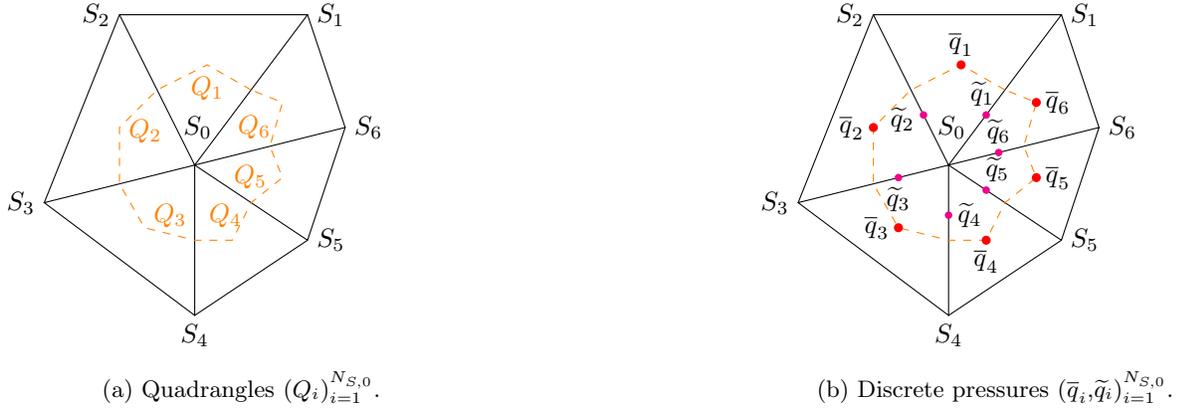
\end{center}
\begin{center}
	\begin{figure}[h]
		\begin{subfigure}[b]{0.39\linewidth}
			\begin{tikzpicture}[scale=0.5]
			\draw[-] (0,0) -- (3,4) -- (-2,4) --cycle;
			\draw[-] (-4,-1) -- (-2,4);
			\draw[-] (-4,-1) -- (4,1);
			\draw[-] (4,1) -- (3,4);
			\draw[dashed,orange] (1.5,2)--(1/3,8/3)--(-1,2)--(-2,1)--(-2,-0.5);
			\draw[dashed,orange] (2,0.5)--(7/3,5/3)--(1.5,2);
			\coordinate [label=left:$S_0$] (S) at (0.7,1);
			\coordinate [label=right:$S_2$] (S2) at (3,4);
			\coordinate [label=left:$S_3$] (S3) at (-2,4);
			\coordinate [label=left:$S_4$] (S4) at (-4,-1);
			\coordinate [label=right:$S_1$] (S1) at (4,1);
			%
			\coordinate [label=below:\textcolor{orange}{$Q_2$}] (Q2) at (1/3,7.8/3);
			\coordinate [label=right:\textcolor{orange}{$Q_3$}] (Q3) at (-2,0.9);
			\coordinate [label=right:\textcolor{orange}{$Q_1$}] (Q1) at (0.90,1);
			\coordinate [label=below:$\tilde{F}_4$] (FI4) at (-1,-0.25);
			\coordinate [label=below:$\tilde{F}_1$] (FI1) at (6/4,0.35);
			\draw[dashed,red] (-2,-0.5)--(0,0);
			\draw[dashed,red] (2,0.5)--(0,0);
			\end{tikzpicture}
			\caption{Quadrangles $(Q_i)_{i=1}^{N_{S,0}}$.}
		\end{subfigure}
		\hspace{0.2\linewidth}
		\begin{subfigure}[b]{0.39\linewidth} 
			\begin{tikzpicture}[scale=0.5]
			\draw[-] (0,0) -- (3,4) -- (-2,4) --cycle;
			\draw[-] (-4,-1) -- (-2,4);
			\draw[-] (-4,-1) -- (4,1);
			\draw[-] (4,1) -- (3,4);
			\draw[dashed,orange] (1.5,2)--(1/3,8/3)--(-1,2)--(-2,1)--(-2,-0.5);
			\draw[dashed,orange] (2,0.5)--(7/3,5/3)--(1.5,2);
			\coordinate [label=left:$S_0$] (S) at (0.7,1);
			\coordinate [label=right:$S_2$] (S2) at (3,4);
			\coordinate [label=left:$S_3$] (S3) at (-2,4);
			\coordinate [label=left:$S_4$] (S4) at (-4,-1);
			\coordinate [label=right:$S_1$] (S1) at (4,1);
			%
			\coordinate [label = above:$\overline{q}_2$] (pT2) at (1/3,8/3);
			\fill[fill=red] (1/3,8/3) circle (0.12cm);
			\coordinate [label = left:$\overline{q}_3$] (pT3) at (-2,1);
			\fill[fill=red] (-2,1) circle (0.12cm);
			\coordinate [label = right:$\overline{q}_1$] (pT1) at (7/3,5/3) ;
			\fill[fill=red] (7/3,5/3) circle (0.12cm);
			\coordinate [label = above:$\tilde{q}_2$] (p2) at (0.93,4/3);
			\fill[fill=magenta] (1,4/3) circle (0.1cm);
			\coordinate [label = left:$\tilde{q}_3$] (p3) at (-2/3,1.2);
			\fill[fill=magenta] (-2/3,4/3) circle (0.1cm);
			\coordinate [label = below:$\tilde{q}_4$] (p4) at (-4/3,-1/3);
			\fill[fill=magenta] (-4/3,-1/3) circle (0.1cm);
			\coordinate [label = below:$\tilde{q}_1$] (p1) at (4/3,1/3);
			\fill[fill=magenta] (4/3,1/3) circle (0.1cm);
			
			\end{tikzpicture}
			\caption{Discrete pressures ($\overline{q}_i$,$\tilde{q}_i)_{i=1}^{N_{S,0}}$.}
		\end{subfigure}
		\caption{MPFA Scheme for $j \in \Ical_S^b$ and $N_{S,0}=4$ }
		\label{fig:Qj}
	\end{figure}
\end{center}
\vspace{-0.5cm}
Let $q_h\in Q_h$. We set $q_{h|K_\ell}:=\ov{q}_\ell$.
\\
Consider $S_0\in\Om$ (Fig. \ref{fig:Qi}). Let us build a piecewise affine approximation of $q_h$ on each quadrangle $(Q_i)_{i=1}^{N_{S,0}}$ (see Fig. \ref{fig:Qi}-(a)). We call this approximation $\tilde{q}_h$. We first introduce auxiliary discrete pressure values $(\tilde{q}_i)_{i=1}^{N_{S,0}}$ on the thirds of the inner edges of $\Mcal_0$ (see Fig. \ref{fig:Qi}-(b)). For all $j\in\Ical_{S,i}$, we define $\Gcal_i(q_h):=\grad\tilde{q}_{h|Q_i}$, using an integration by part as it is done in \cite[Sect. 3]{LEPOTIER05bis} and \cite[Sect. 1.1.1]{lepotierhdr}: 
$$ |Q_i| \Gcal_i= \int_{Q_i} \Gcal_i(q_h) = \int_{\partial Q_i} \tilde{q}_h  \nvec_{\partial Q_i}=\, \tilde{q}_i \frac{\Scal_{i,i}}{d} + \tilde{q}_{i+1}\,\frac{\Scal_{i+1,i}}{d}  + \ov{q}_i(-\frac{\Scal_{i,i}}{d}-\frac{\Scal_{i+1,i}}{d}). $$

Hence, noticing that $|Q_i|= \frac{|T_i|}{d+1}$, we have:
\begin{equation}\label{eq:gradientlocal1}
\Gcal_i(q_h)=\ds\frac{1}{|Q_i|}\left(\,(\tilde{q}_i-\ov{q}_i)\,\frac{\Scal_{i,i}}{d}+(\tilde{q}_{i+1}-\ov{q}_i)\,\frac{\Scal_{i+1,i}}{d}\right) =\ds\frac{d+1}{d\,|T_i|}\left(\,\tilde{q}_i\,\Scal_{i,i}+\tilde{q}_{i+1}\,\Scal_{i+1,i}+\ov{q}_i\,\Scal_{0,i}\,\right).
\end{equation}
In order to preserve the flux across the inner edges of $\Mcal_0$, we write that:
\begin{equation}\label{eq:continuitefacesnormales}
\forall i\in\Ical_{S,0},\quad \Gcal_i(q_h)\cdot\Scal_{i+1,i}+\Gcal_{i+1}(q_h)\cdot\Scal_{i+1,i+1}=0.
\end{equation}
These $N_{S,0}$ equations with $N_{S,0}$ unknowns (the auxiliary discrete pressure values $(\tilde{q}_i)_{i=1}^{N_{S,0}}$) lead to a well posed linear system. Thus, we can evaluate the auxiliary discrete pressure values $(\tilde{q}_i)_{i=1}^{N_{S,0}}$ with the data $(\ov{q}_i)_{i=1}^{N_{S,0}}$. Therefore, we can explicitly express the pressure gradients $(\Gcal_i(q_h))_{i=1}^{N_{S,0}}$ \eqref{eq:gradientlocal1}.
\\
Consider now $S_0\in\pa\Om$ (see Fig. \ref{fig:Qj}).  According to \cite[proof of Prop. IV.3.7]{BoyerFabrie12}, if $\fvec\in\bH^1(\Om)$, the solution $(\uvec,p)$ to Problem \eqref{eq:Stokes} is such that:
\begin{equation}\label{eq:Stokes-FVh-P0}
\grad p \cdot \nvec_{|\partial \Omega}=\fvec\cdot \nvec_{|\partial \Omega},
\end{equation}
where $\nvec_{|\partial \Omega}$ is the unit outward normal vector at $\partial \Omega$. 

In our numerical experiments, we explicit the auxiliary discrete pressure values located on $\pa\Om$ (ie $\tilde{q}_1$ and $\tilde{q}_4$ on Fig. \ref{fig:Qj}-(b)) by imposing that for all $i\in\Ical_{S,0}$ such that $F_i\in\pa\Om$:
\begin{equation}
\label{eq:grap.n=f.n}
\int_{\tilde{F}_i}  \Gcal_i(q_h)\cdot\nvec_{|\tilde{F}_i}= \int_{\tilde{F}_i} \fvec\cdot \nvec_{|\tilde{F}_i}.
\end{equation}
Again, the auxiliary discrete pressure values solve a well posed linear system. They can be written with the data $(\ov{q}_i)_{i=1}^{N_{S,0}}$ and we can explicitly express $\Gcal_i(q_h)$.
\\

For $i\in\Ical_S$, we let $(Q_{i,j})_j\in\Ical_{S,i}$ be the set of quadrangles built around $S_i$, and we call $\mathcal{Q}_h$ the mesh of all the quadrangles $\mathcal{Q}_h:=(\,(Q_{i,j})_{j\in\Ical_{S,i}})_{i\in\Ical_S}$. Let $q_h\in Q_h$. Let $i\in\Ical_S$. In the macro-element $\Mcal_i$, we call $\Gcal_{i,j}(q_h)$ the local reconstructed gradient of $q_h$. We now define the MPFA gradient reconstruction as the operator $\Gcal_h$: 

\begin{equation}
\label{eq:gradientMPFA}
\Gcal_h:\left\{
\begin{array}{rcl}Q_h&\rightarrow&\bP^0(\mathcal{Q}_h)\\
q_h&\mapsto&\Gcal_h(q_h)
\end{array}\right.
\,|\quad  \forall  i\in\Ical_S, \,  \forall j \in \Ical_{S,i},\quad
\Gcal_h(q_h)_{|Q_{i,j}} = \Gcal_{i,j}(q_{h|\Mcal_i}).
\end{equation}
If the data $\fvec$ is of low regularity, one can enhance the space of discrete pressures, adding the auxiliary unknowns on the boundary as degrees of freedom.
\\

\begin{prop}
	With triangles, and $p\in C^2(\Om)$, the fluxes of the symmetric MPFA scheme are consistently approximated. Also by choosing the auxiliary pressures unknowns at the thirds of the edges, the gradient is approximated exactly for affine functions. Also, the symmetric MPFA scheme is  consistent, coercive and convergent. 
\end{prop} 
The proof of this proposition can be found in \cite[Prop. 2, Prop. 3 ]{LEPOTIER05bis} and \cite[Theorem 3.2]{LiShYo09}.\\

Let us express our discrete Stokes problem. Let $g_h(\cdot,\cdot)$ be the following bilinear form:
\begin{equation}\label{eq:bilin-form-gh}
g_h:\left\{
\begin{array}{rcl}
\bX_{h}\times Q_h&\rightarrow&\R\\
(\vvec_h,q_h)&\mapsto&(\Gcal_h(q_h),\vvec_h)_{\bL^2(\Om)}
\end{array}
\right..
\end{equation}

The discretization of \eqref{eq:Stokes} using the MPFA scheme to discretize the pressure gradient reads:
\begin{equation}
\label{eq:MPFA-VF}
\mbox{Find }(\uvec,p)\in\bX_{0,h}\times Q_h\,|\quad \left\{
\begin{array}{rcll}
a_{\nu,h}(\uvec_h, \vvec_h)+g_h(\vvec_h, p_h)&=&(\fvec, \vvec_h)_{\bL^2(\Om)} & \forall \vvec_h \in  \textbf{X}_h\\
b_h(\uvec_h, q_h)&=&0 &\forall q_h \in Q_h
\end{array}
\right.,
\end{equation}
where the bilinear forms $a_{\nu,h}(\cdot,\cdot)$ and $b_h(\cdot,\cdot)$ are defined by \eqref{eq:DiscBilinForms}. Notice that the linear system related to variational formulation \eqref{eq:MPFA-VF} is not symmetric.

\section{Numerical Results on the Stokes Problem}\label{sec:numstokes}
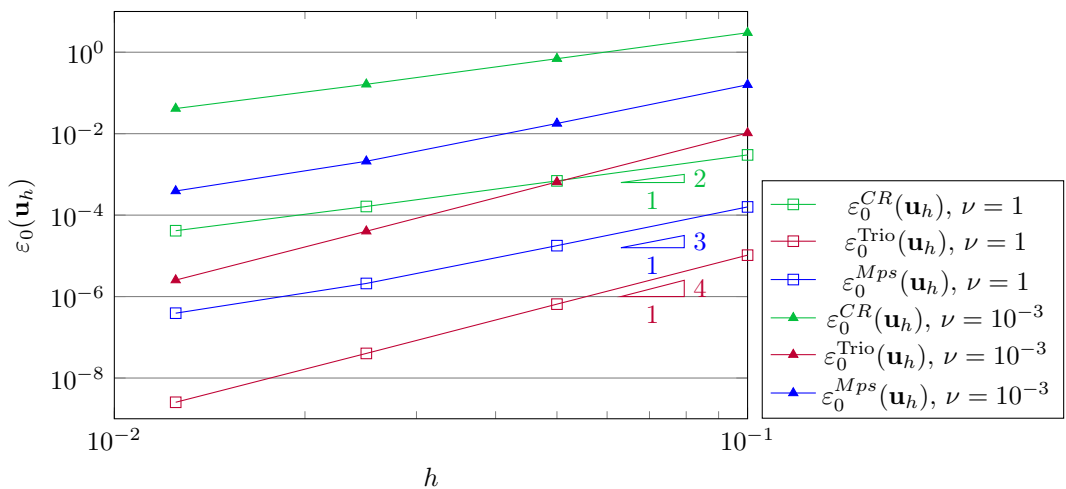
\begin{figure}[!ht]
	\centering
	\begin{tikzpicture}
	\begin{axis}[
	height = 7cm, 
	width = 10cm,  
	xlabel = {$h$},
	ylabel = {$\eps_0(\uvec_h)$},
	ymajorgrids=true,
	major grid style={black!50},
	xmode=log, ymode=log,
	xmin=1e-2, xmax=1e-1,
	ymin=1e-9, ymax=1e1,
	yticklabel style={
		/pgf/number format/fixed,
		/pgf/number format/precision=0},
	scaled y ticks=false,
	xticklabel style={
		/pgf/number format/fixed,
		/pgf/number format/precision=0},
	scaled y ticks=false,  
	legend style={at={(1.5,0)},anchor=south east},
	legend columns=1
	], 
	\addplot[color=mygreen,mark=square] table [x=h, y=eUCRP0nu1] {Plots/Ctu0psin.txt};
	\addlegendentry{$\eps_0^{CR}(\uvec_h) ,\, \nu=1$ }
	\addplot[color=myred,mark=square] table [x=h, y=eUCRP0P1nu1] {Plots/Ctu0psin.txt};
	\addlegendentry{$\eps_0^{\textrm{Trio}}(\uvec_h),\, \nu=1 $}
	\addplot[color=myblue,mark=square] table [x=h, y=eUCRMPSnu1] {Plots/Ctu0psin.txt};
	\addlegendentry{$\eps_0^{Mps}(\uvec_h) ,\, \nu=1$}
	\addplot[color=mygreen,mark=triangle*] table [x=h, y=eUCRP0nu0001] {Plots/Ctu0psin.txt};
	\addlegendentry{$\eps_0^{CR}(\uvec_h),\, \nu=10^{-3} $ }
	\addplot[color=myred,mark=triangle*] table [x=h, y=eUCRP0P1nu0001] {Plots/Ctu0psin.txt};
	\addlegendentry{$\eps_0^{\textrm{Trio}}(\uvec_h) ,\, \nu=10^{-3}$}
	\addplot[color=myblue,mark=triangle*] table [x=h, y=eUCRMPSnu0001] {Plots/Ctu0psin.txt};
	\addlegendentry{$\eps_0^{Mps}(\uvec_h) ,\, \nu=10^{-3}$}
	\logLogSlopeTriangle{0.9}{0.1}{0.42}{3}{myblue};
	\logLogSlopeTriangle{0.9}{0.1}{0.58}{2}{mygreen};
	\logLogSlopeTriangle{0.9}{0.1}{0.30}{4}{myred};
	\end{axis}       
	\end{tikzpicture}
	\caption{$\eps_0^X(\uvec_h)$ for $(\uvec,p) =(\mathbf{0},\sin(2 \pi x )\,\sin( 2 \pi y)\,)$}
	\label{fig:NoFlowu1}
\end{figure}

\begin{figure}[!ht]
	\centering
	\begin{tikzpicture}
	\begin{axis}[
	height = 7cm, 
	width = 10cm,  
	xlabel = {$h$},
	ylabel = {$\eps_0(p_h)$},
	ymajorgrids=true,
	major grid style={black!50},
	xmode=log, ymode=log,
	xmin=1e-2, xmax=1e-1,
	ymin=1e-4, ymax=1e-0,
	yticklabel style={
		/pgf/number format/fixed,
		/pgf/number format/precision=0},
	scaled y ticks=false,
	xticklabel style={
		/pgf/number format/fixed,
		/pgf/number format/precision=0},
	scaled y ticks=false,  
	legend style={at={(1.5,0)},anchor=south east},
	legend columns=1
	],
	\addplot[color=mygreen,mark=square] table [x=h, y=ePCRP0nu1] {Plots/Ctu0psin.txt};
		\addlegendentry{$\eps_0^{CR}(p_h) ,\, \nu=1$ }
		\addplot[color=myred,mark=square] table [x=h, y=ePCRP0P1nu1] {Plots/Ctu0psin.txt};
		\addlegendentry{$\eps_0^{\textrm{Trio}}(p_h),\, \nu=1 $}
		\addplot[color=myblue,mark=square] table [x=h, y=ePCRMPSnu1] {Plots/Ctu0psin.txt};
		\addlegendentry{$\eps_0^{Mps}(p_h) ,\, \nu=1$}
	\addplot[color=mygreen,mark=triangle*] table [x=h, y=ePCRP0nu0001] {Plots/Ctu0psin.txt};
		\addlegendentry{$\eps_0^{CR}(p_h),\, \nu=10^{-3} $ }
	\addplot[color=myred,mark=triangle*] table [x=h, y=ePCRP0P1nu0001] {Plots/Ctu0psin.txt};
		\addlegendentry{$\eps_0^{\textrm{Trio}}(p_h) ,\, \nu=10^{-3}$}
	\addplot[color=myblue,mark=triangle*] table [x=h, y=ePCRMPSnu0001] {Plots/Ctu0psin.txt};
		\addlegendentry{$\eps_0^{Mps}(p_h) ,\, \nu=10^{-3}$}
		\logLogSlopeTriangle{0.9}{0.1}{0.43}{2}{myred};
		\logLogSlopeTriangle{0.9}{0.1}{0.70}{1}{mygreen};
		\logLogSlopeTriangle{0.9}{0.1}{0.32}{2}{myblue};
	\end{axis}       
	\end{tikzpicture}
	\caption{$\eps_0(p_h)$ for $(\uvec,p) =(\mathbf{0},\sin(2 \pi x )\,\sin( 2 \pi y)\,)$}
	\label{fig:NoFlowp1}
\end{figure}
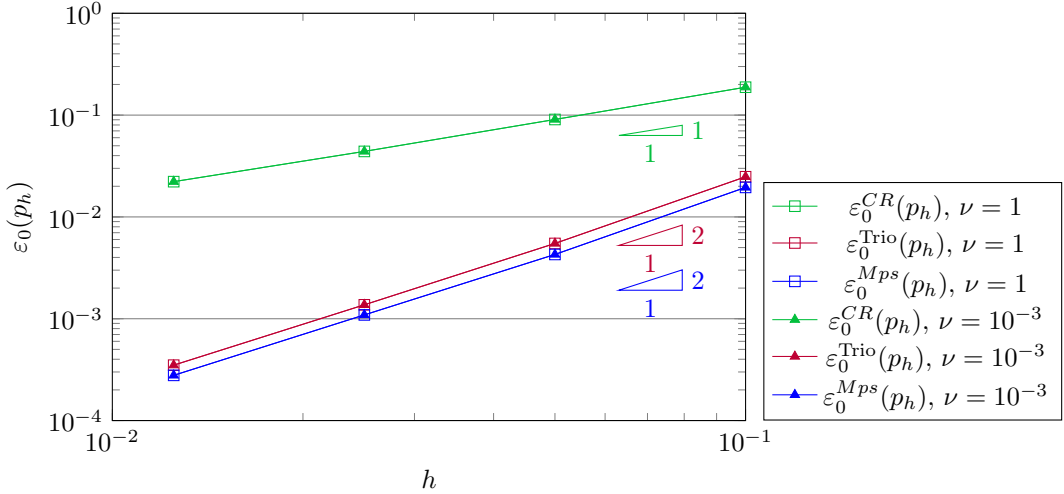

In this section, we give some $2D$ numerical results which compare the $\bP^1_{nc}-P^0_{Mps}$ scheme to the $\bP^1_{nc}-P^0$ and $\bP^1_{nc}-(P^0+P^1)$ schemes. \\
Consider Problem \eqref{eq:StokesVF} with prescribed solution such that: $ (\uvec,p)=(\bm{0},\varphi)$. \\
When $\varphi$ is some affine function, then both  $\bP^1_{nc}-(P^0+P^1)$ and  $\bP^1_{nc}-P^0_{Mps}$ schemes give exactly $\uvec_h= \bm{0}$.\\
When $\varphi$ is some quadratic function, then $\bP^1_{nc}-(P^0+P^1)$ scheme gives exactly $\uvec_h= \bm{0}$, as a consequence of Property \ref{thm:P1nc-P0P1}. \\
In what follows, we set $\Om=(0,1)^2$. We denote the $L^2(\Om)$ errors estimates of the discrete velocity and pressure by:
$$\eps_0^X(\uvec_h) := \left\{ 
\begin{array}{rl}
\ds \|\uvec_h\|_{\bL^2(\Om)}&\mbox{ if }\uvec=0\\ \\
 \ds \frac{\|\uvec_h-\uvec\|_{\bL^2(\Om)}}{ \|\uvec \|_{\bL^2(\Om)}} &\mbox{ otherwise}\end{array}\right.  \quad  \quad   \textrm{and} \quad \quad    \eps_0^X(p_h) := \frac{\| p_h-p\|_{L^2(\Om)}}{\|p\|_{L^2(\Om)}}, $$
 where: $X=CR$ (resp. $X=Trio$ and $X=Mps$) refers to the solution computed with the $\bP^1_{nc}-P^0$ (resp. $\bP^1_{nc}-(P^0+P^1)$ and $\bP^1_{nc}-P^0_{Mps}$) scheme. \\
 We first consider Problem \eqref{eq:StokesVF} with prescribed solution $(\uvec,p) =(\mathbf{0},\sin(2 \pi x )\,\sin( 2 \pi y)\,)$. On Fig. \ref{fig:NoFlowu1} (resp. \ref{fig:NoFlowp1}), we plot $\eps_0^X(\uvec_h)$ (resp. $\eps_0^X(p_h)$) against the meshstep $h$ in the logarithmic scale, for $\nu=1$ and $\nu=10^{-3}$.
 
  We notice that $\eps_0^{X}(\uvec_h)\propto\nu^{-1}$ for the three schemes. Concerning the $\bP^1_{nc}-P^0_{Mps}$ scheme, we first remark that $\eps_0^{Mps}(\uvec_h)$ gives intermediate results between $\eps_0^{CR}(\uvec_h)$ and $\eps_0^{Trio}(\uvec_h)$ (see Fig. \ref{fig:NoFlowu1}). Second, we notice that $\eps_0^{Mps}(p_h)\approx\eps_0^{Trio}(p_h)$. Finally, we notice that the $\bP^1_{nc}-P^0_{Mps}$ scheme returns a convergence rate of order $3$ for $\eps_0^{Mps}(\uvec_h)$ and $2$ for $\eps_0^{Mps}(p_h)$.

\newpage

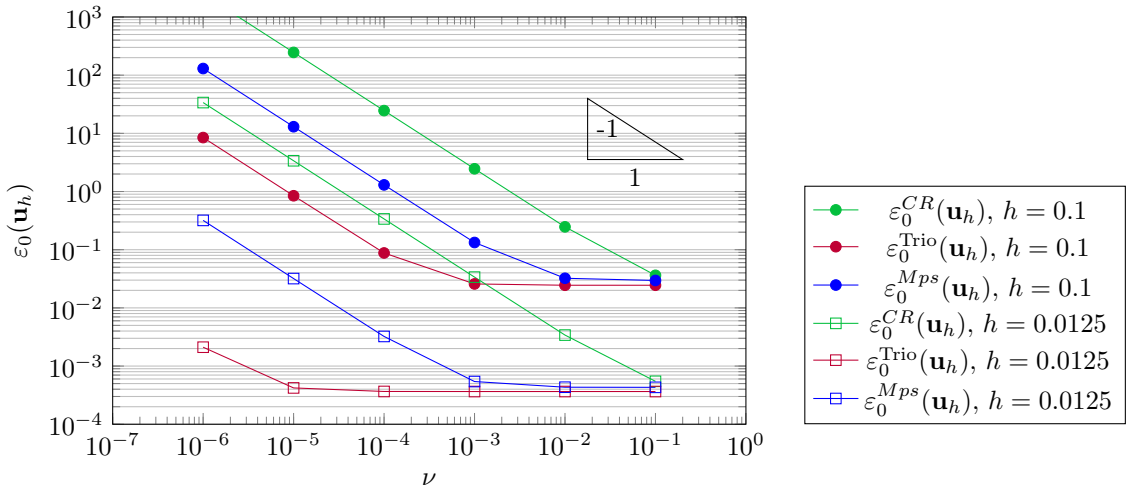
\begin{figure}[!ht]
	\centering
	\begin{tikzpicture}
	\begin{axis}[
	height = 7cm, 
	width = 10cm,  
	xlabel = {$\nu$},
	ylabel = {$\eps_0(\uvec_h)$},
	yminorgrids=true,
	major grid style={black!50},
	xmode=log, ymode=log,
	xmin=1e-7, xmax=1e-0,
	ymin=1e-4, ymax=1e03,
	yticklabel style={
		/pgf/number format/fixed,
		/pgf/number format/precision=0},
	scaled y ticks=false,
	xticklabel style={
		/pgf/number format/fixed,
		/pgf/number format/precision=0},
	scaled y ticks=false,  
	legend style={at={(1.6,0)},anchor=south east},
	legend columns=1
	],
	\addplot[color=mygreen,mark=*] table [x=visco, y=eUCRP0] {Plots/visco/visco_Stokes_num30_h1TrioCFDconvnu.txt};
	\addlegendentry{$\eps_0^{CR}(\uvec_h),\, h=0.1 $ }
	\addplot[color=myred,mark=*] table [x=visco, y=eUCRP0P1] {Plots/visco/visco_Stokes_num30_h1TrioCFDconvnu.txt};
	\addlegendentry{$\eps_0^{\textrm{Trio}}(\uvec_h),\, h=0.1 $}
	\addplot[color=myblue,mark=*] table [x=visco, y=eUCRMPS] {Plots/visco/visco_Stokes_num30_h1TrioCFDconvnu.txt};
	\addlegendentry{$\eps_0^{Mps}(\uvec_h),\, h=0.1 $}	
	\addplot[color=mygreen,mark=square] table [x=visco, y=eUCRP0] {Plots/visco/visco_Stokes_num30_h4TrioCFDconvnu.txt};
	\addlegendentry{$\eps_0^{CR}(\uvec_h),\, h=0.0125 $ }
	\addplot[color=myred,mark=square] table [x=visco, y=eUCRP0P1] {Plots/visco/visco_Stokes_num30_h4TrioCFDconvnu.txt};
	\addlegendentry{$\eps_0^{\textrm{Trio}}(\uvec_h), \, h=0.0125 $}
	
	\addplot[color=myblue,mark=square] table [x=visco, y=eUCRMPS] {Plots/visco/visco_Stokes_num30_h4TrioCFDconvnu.txt};
	\addlegendentry{$\eps_0^{Mps}(\uvec_h),\, h=0.0125 $}
	\logLogSlopeTriangle{0.75}{-0.15}{0.65}{-1}{black};
	\end{axis}       
	\end{tikzpicture}
	\caption{$\eps_0(\uvec_h)$ for $\uvec$ and $p$ sinusoidal functions against viscosity with different mesh sizes}
	\label{fig:testvisco_vit}
\end{figure}

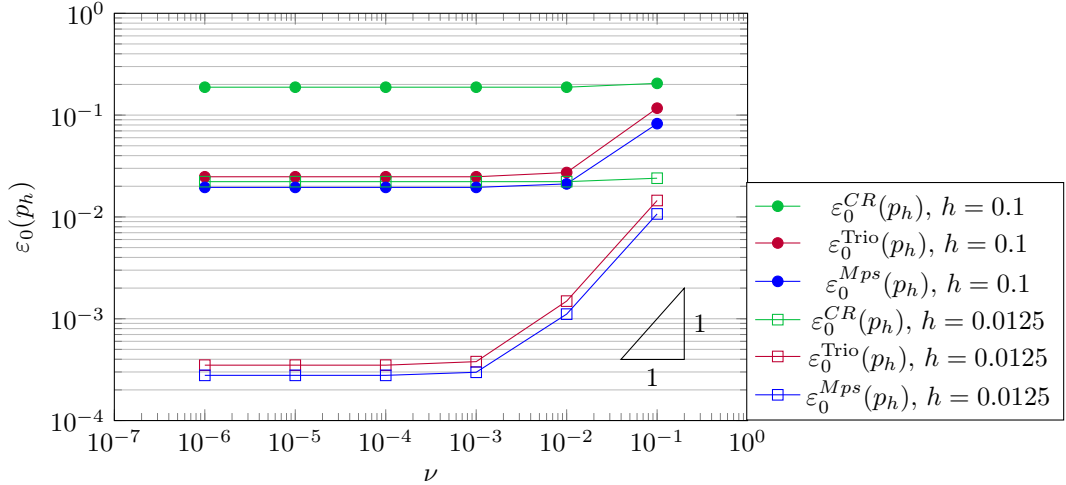
\begin{figure}[!ht]
	\centering
	\begin{tikzpicture}
	\begin{axis}[
	height = 7cm, 
	width = 10cm,  
	xlabel = {$\nu$},
	ylabel = {$\eps_0(p_h)$},
	yminorgrids=true,
	major grid style={black!50},
	xmode=log, ymode=log,
	xmin=1e-7, xmax=1e-0,
	ymin=1e-4, ymax=1e0,
	yticklabel style={
		/pgf/number format/fixed,
		/pgf/number format/precision=0},
	scaled y ticks=false,
	xticklabel style={
		/pgf/number format/fixed,
		/pgf/number format/precision=0},
	scaled y ticks=false,  
	legend style={at={(1.5,0)},anchor=south east},
	legend columns=1
	],
	\addplot[color=mygreen,mark=*] table [x=visco, y=ePCRP0] {Plots/visco/visco_Stokes_num30_h1TrioCFDconvnu.txt};
	\addlegendentry{$\eps_0^{CR}(p_h),\, h=0.1  $ }
	\addplot[color=myred,mark=*] table [x=visco, y=ePCRP0P1] {Plots/visco/visco_Stokes_num30_h1TrioCFDconvnu.txt};
	\addlegendentry{$\eps_0^{\textrm{Trio}}(p_h),\, h=0.1  $}
	\addplot[color=myblue,mark=*] table [x=visco, y=ePCRMPS] {Plots/visco/visco_Stokes_num30_h1TrioCFDconvnu.txt};
	\addlegendentry{$\eps_0^{Mps}(p_h),\, h=0.1  $}
	\addplot[color=mygreen,mark=square] table [x=visco, y=ePCRP0] {Plots/visco/visco_Stokes_num30_h4TrioCFDconvnu.txt};
	\addlegendentry{$\eps_0^{CR}(p_h),\, h=0.0125  $ }
	\addplot[color=myred,mark=square] table [x=visco, y=ePCRP0P1] {Plots/visco/visco_Stokes_num30_h4TrioCFDconvnu.txt};
	\addlegendentry{$\eps_0^{\textrm{Trio}}(p_h),\, h=0.0125  $}
	\addplot[color=myblue,mark=square] table [x=visco, y=ePCRMPS] {Plots/visco/visco_Stokes_num30_h4TrioCFDconvnu.txt};
	\addlegendentry{$\eps_0^{Mps}(p_h),\, h=0.0125  $}
	\logLogSlopeTriangle{0.9}{0.1}{0.15}{1}{black};
	\end{axis}       
	\end{tikzpicture}
	\caption{$\eps_0(p_h)$ for $\uvec$ and $p$ sinusoidal functions against viscosity with different mesh sizes}
	\label{fig:testvisco_press}
\end{figure}

We notice that, compared to the $P^1_{NC}-P^0$ scheme, the errors are greatly reduced by $\bP^1_{nc}-(P^0+P^1)$ and $\bP^1_{nc}-P^0_{Mps}$ schemes. These schemes allow to attenuate the amplitude of spurious velocities and hence provide a better simulation. This is illustrated by the resolution of \eqref{eq:Stokes} with $(\uvec,p)$ defined by \eqref{Test case 4 }. In this case, as $\uvec$ is not an affine function, the three schemes return a convergence rate of order $2$ for $\eps_0^{Mps}(\uvec_h)$ and $1$ for $\eps_0^{Mps}(p_h)$. The errors resulted for $h=0.1$ and $h=0.0125$ are plotted against viscosity in Figures \ref{fig:testvisco_vit} and \ref{fig:testvisco_press}. In these plots, we see that the $\bP^1_{nc}-P^0_{Mps}$ scheme gives intermediate results. Also, we notice that the spurious velocities errors become overriding when: 
\begin{itemize}
	\item $\nu\leq 10^0$    with $h=0.1$  and $\nu\leq 10^{-3}$ with $h=0.0125$ for the $P^1_{NC}-P^0$. 
	\item $\nu\leq 10^{-2}$ with $h=0.1$  and $\nu\leq 10^{-3}$ with $h=0.0125$ for the $\bP^1_{nc}-P^0_{Mps}$. 
	\item $\nu\leq 10^{-3}$ with $h=0.1$  and $\nu\leq 10^{-5}$ with $h=0.0125$ for the $\bP^1_{nc}-(P^0+P^1)$. 
\end{itemize}

The tipping viscosity point, where the spurious velocities errors become dominant, depends on the velocity error generated by the gradient approximation and therefore the mesh size. As these schemes converge with different orders when $\uvec=0$, it can be seen that decreasing the mesh size reduces the viscosity at which this point is reached more or less depending on the order. 

\begin{equation} 
\label{Test case 4 }
 (\uvec,p)= \begin{pmatrix} \begin{matrix}
(\cos(2 \pi x)-1) \, \sin(2 \pi y) \\ 
-(\cos(2 \pi y)-1)\, \sin(2 \pi x)
\end{matrix} \, 
,  \sin(2 \pi x ) \sin( 2 \pi y)
\end{pmatrix}
\end{equation}

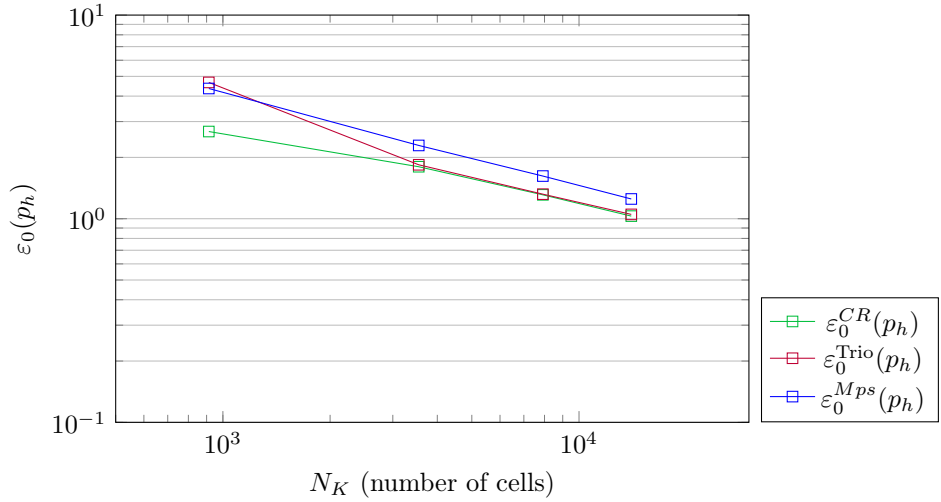
\begin{figure}[!ht]
	\centering
	\begin{tikzpicture}
	\begin{axis}[
	height = 7cm, 
	width = 10cm,  
	xlabel = {$N_K$ (number of cells)},
	ylabel = {$\eps_0(p_h)$},
	yminorgrids=true,
	major grid style={black!50},
	xmode=log, ymode=log,
	xmin=5e2, xmax=3e4,
	ymin=1e-1, ymax=1e1,
	yticklabel style={
		/pgf/number format/fixed,
		/pgf/number format/precision=0},
	scaled y ticks=false,
	xticklabel style={
		/pgf/number format/fixed,
		/pgf/number format/precision=0},
	scaled y ticks=false,  
	legend style={at={(1.3,0)},anchor=south east},
	legend columns=1
	],
	\addplot[color=mygreen,mark=square] table [x=NCR, y=ePCRP0] {Plots/anisotropic/Stokes_meshtype8_num30ef0_2_7.txt};
	\addlegendentry{$\eps_0^{CR}(p_h) $ }
	\addplot[color=myred,mark=square] table [x=NCR, y=ePCRP0P1] {Plots/anisotropic/Stokes_meshtype8_num30ef0_2_7.txt};
	\addlegendentry{$\eps_0^{\textrm{Trio}}(p_h)$}
	\addplot[color=myblue,mark=square] table [x=NCR, y=ePCRMPS] {Plots/anisotropic/Stokes_meshtype8_num30ef0_2_7.txt};
	\addlegendentry{$\eps_0^{Mps}(p_h) $}
	\end{axis}       
	\end{tikzpicture}
	\caption{$\eps_0(p_h)$ for $\uvec$ and $p$ sinusoidal functions against viscosity with different kershaw meshes and $\nu =1$.}
	\label{fig:kershawp}
\end{figure}

\newpage
We are also interested in the sensitivity to the mesh deformation. Indeed, nowadays, mesh refinement techniques based on a posteriori error estimators or industrial constraints can generate anisotropic meshes. In this subsection, we show that the three schemes have the same behaviour with respect to the regularity of the mesh.
To illustrate this property, we propose to use the Kershaw meshes presented in the benchmark \cite{HeHu08} ( see Fig. \ref{fig:kershaw}) with $(\uvec,p)$ in \eqref{Test case 4 }. As the mesh is composed of quadrilaterals, we cut them with along one of the diagonal, which allows us to remain within reasonable convergence assumptions. The mesh is represented in Fig. \ref{fig:kershaw} and we plot the results in Fig \ref{fig:kershawu} and \ref{fig:kershawp}. We can see that the schemes have a convergence rate of order $2$ for $\eps_0^{X}(\uvec_h)$ and $1$ for $\eps_0^{X}(p_h)$.

\begin{figure}[!ht]
	\centering
	\includegraphics[width=7.2cm]{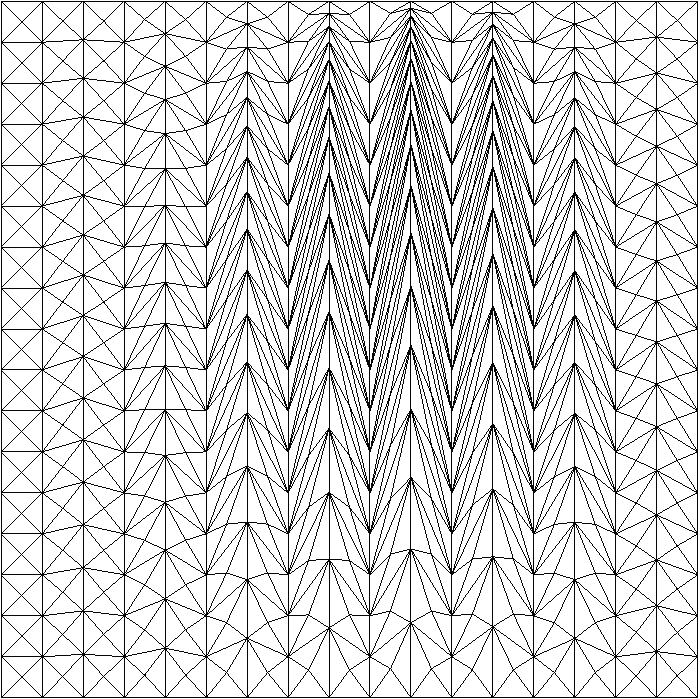}
	\caption{Kershaw mesh.}
	\label{fig:kershaw}
\end{figure} 
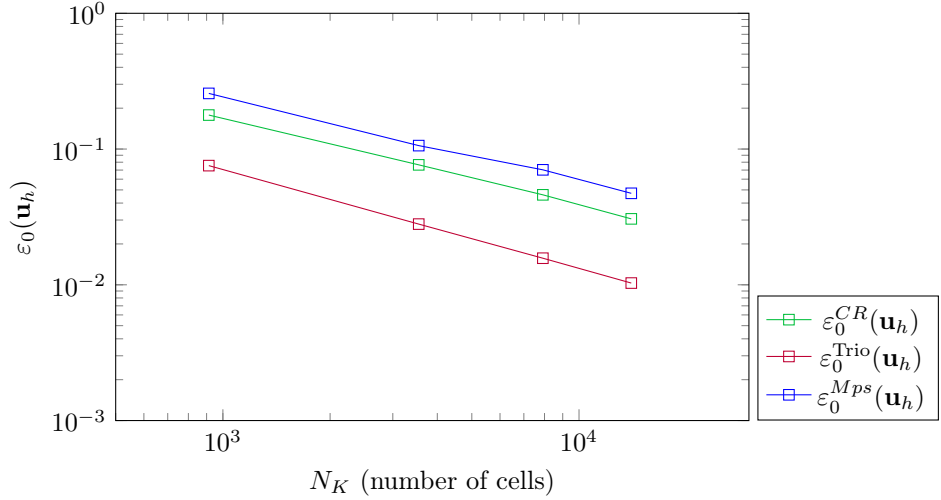
\begin{figure}[!ht]
	\centering
	\begin{tikzpicture}
	\begin{axis}[
	height = 7cm, 
	width = 10cm,  
	xlabel = {$N_K$ (number of cells)},
	ylabel = {$\eps_0(\uvec_h)$},
	major grid style={black!50},
	xmode=log, ymode=log,
	xmin=5e2, xmax=3e4,
	ymin=1e-3, ymax=1e0,
	yticklabel style={
		/pgf/number format/fixed,
		/pgf/number format/precision=0},
	scaled y ticks=false,
	xticklabel style={
		/pgf/number format/fixed,
		/pgf/number format/precision=0},
	scaled y ticks=false,  
	legend style={at={(1.3,0)},anchor=south east},
	legend columns=1
	], 
	\addplot[color=mygreen,mark=square] table [x=NCR, y=eUCRP0] {Plots/anisotropic/Stokes_meshtype8_num30ef0_2_7.txt};
	\addlegendentry{$\eps_0^{CR}(\uvec_h) $ }
	\addplot[color=myred,mark=square] table [x=NCR, y=eUCRP0P1] {Plots/anisotropic/Stokes_meshtype8_num30ef0_2_7.txt};
	\addlegendentry{$\eps_0^{\textrm{Trio}}(\uvec_h) $}
	\addplot[color=myblue,mark=square] table [x=NCR, y=eUCRMPS] {Plots/anisotropic/Stokes_meshtype8_num30ef0_2_7.txt};
	\addlegendentry{$\eps_0^{Mps}(\uvec_h) $}
	\end{axis}       
	\end{tikzpicture}
	\caption{$\eps_0(\uvec_h)$ for $\uvec$ and $p$ sinusoidal functions against viscosity with different kershaw meshes and $\nu =1$. }
	\label{fig:kershawu}
\end{figure}

\newpage
\newpage

\section{Numerical Results on the Navier-Stokes Problem}\label{sec:numNstokes}
 We choose the convection scheme initially presented in \cite[Eq. 2.8]{GaHeLa10}. This choice is motivated by the result of the Benchmark \cite{CHEY21} where the scheme presents good convergence and stability results  without increasing the stencil of the scheme. To resolve efficiently the Navier-Stokes equations, we use a prediction correction time-scheme  \cite{Chorin67} \cite{Teman68},  which consists in calculating a predicted velocity (with non-zero divergence), then solving the pressure at the next time and correcting the velocity to ensure divergence-free flow. With this approach, the velocity and pressure resolutions are decoupled.

	 As the approach is not classical, it is interesting to check the convergence of the scheme on Navier-Stokes. This leads to study the Green-Taylor vortex solution which is a well-known analytical solution to \eqref{eq:NS}.

Let first introduce the space-discretization of the Navier-Stokes equations \eqref{eq:NS} : 

\begin{equation}
\left\{
\begin{aligned}
 M \partial_t U +   \nu KU +L(U)+GP &=F\\
DU&=0
\end{aligned}
\right.
\end{equation}

where $U$, $P$ contains the velocity and pressure unknowns and $F$ is the right hand side. The matrices $M$ and $K$ are respectively the mass and stiffness matrices. Also, the matrices $G$ and $D$ represent the gradient and divergence operators. Finally, the matrix $L(U)$ is associated with the convection term and $\partial_t U$ is the time derivative of $U$.

We first present the prediction-correction time scheme : 

 \begin{enumerate}
 	\item Prediction :
 	\begin{equation}
 	\label{eq:pred_corr1}
 	M \frac{U^*-U^n}{\delta t} +\nu KU^* + L(U^n)U^n+GP^n =F^n \quad \mbox{ and } \quad DU^*\neq 0
 	\end{equation}
 	\item Pressure solver:
 	\begin{equation}
 	\label{eq:pred_corr2}
 	\delta t (D \tilde{M}^{-1} G) \delta P = D U^*
 	\end{equation}
 	with $\delta P = P^{n+1}-P^n$
 	\item Correction:
 	\begin{equation}
 	\label{eq:pred_corr3}
 	U^{n+1} = U^* +\delta t \tilde{M}^{-1} G \delta P
 	\end{equation}
 \end{enumerate}
with $\tilde{M}=M+\delta t \, \nu  K$. The CFL of the global system is then: 
\begin{equation}
	\label{eq:CFL_ORDER1}
	\delta t < C h 
\end{equation}

\begin{rmq}
	In Equation \eqref{eq:pred_corr2}, we can approximate $\tilde{M} \simeq M $. This leads to a linear system which is faster to solve but less accurate. 
\end{rmq} 

\begin{rmq}
	As we did in section 6, to determine the auxiliary pressures on the boundary of the MPFA scheme, we impose a condition for the pressure gradient at the edge. On each half-edge $\tilde{F}_i$ related to the vertex $S_0$ and edge $F_i \in \partial \Om$:
	\begin{equation}
	\label{eq:grap.n=f.n_for_NS}
	\int_{\tilde{F}_i}  \Gcal_i\cdot\nvec_{|\tilde{F}_i}= \int_{\tilde{F}_i} \big(\fvec + (\uvec_h^{n+1}-\uvec_h^{n})/\delta t + (\uvec_h^n \cdot \grad) \uvec_h^n \big)\cdot \nvec_{|\tilde{F}_i}.
	\end{equation}
\end{rmq}

Let $\Om=(0,1)^2$. We prescribe the solution to Equation \eqref{eq:NS} with $\fvec=\bm{0}$, to be:
 
 \begin{equation}
 \label{Test case Green-Taylor }
 \left\{
 \begin{array}{rcll}
 u_x &=& - \, \cos\big(2 \pi (x+\frac{1}{2})\big) \,\sin\big(2 \pi (y+\frac{1}{2})\, \big)exp(-8\pi^2t) \vspace{0.05cm}\\
 u_y &=&  \, \sin\big(2 \pi (x+\frac{1}{2})\big)\, \cos \big(2 \pi (y+\frac{1}{2})\, \big)exp(-8\pi^2t) \vspace{0.05cm} \\
 p &=&   -\frac{1}{4} \, \cos(4 \pi (x+\frac{1}{4}) ) + \cos( 4 \pi (y+\frac{1}{2}))\, exp(-16\pi^2t) \vspace{0.05cm}
 \end{array}
 \right.
 \end{equation}
 
 We set $t_{max}=0.01 $ the final time of the simulation. The time step is chosen with respect to the CFL \eqref{eq:CFL_ORDER1} with $C=\frac{1}{2}$. The errors in velocity and pressure at the final time are plotted in Figures \ref{fig:NSgreentaylorU} \ref{fig:NSgreentaylorP} against the mesh step. We can see that the three schemes converge with order $2$ for $\eps_0^{Mps}(\uvec_h)$ and $1$ for $\eps_0^{Mps}(p_h)$ as expected.

\begin{figure}[!ht]
	\centering
	\begin{tikzpicture}
	\begin{axis}[
	height = 7cm, 
	width = 10cm,  
	xlabel = {$h$},
	ylabel = {$\eps_0(\uvec_h)$},
	ymajorgrids=true,
	major grid style={black!50},
	xmode=log, ymode=log,
	xmin=1e-2, xmax=1e-1,
	ymin=1e-4, ymax=1e-1,
	yticklabel style={
		/pgf/number format/fixed,
		/pgf/number format/precision=0},
	scaled y ticks=false,
	xticklabel style={
		/pgf/number format/fixed,
		/pgf/number format/precision=0},
	scaled y ticks=false,  
	legend style={at={(1.3,0)},anchor=south east},
	legend columns=1
	], 
	\addplot[color=mygreen,mark=triangle*] table [x=h, y=eUCRP0] {Plots/NavierStokes_num100ef0_2_7.txt};
	\addlegendentry{$\eps_0^{CR}(\uvec_h) $ }
	\addplot[color=myred,mark=triangle*] table [x=h, y=eUCRP0P1] {Plots/NavierStokes_num100ef0_2_7.txt};
	\addlegendentry{$\eps_0^{\textrm{Trio}}(\uvec_h)$}
	\addplot[color=myblue,mark=triangle*] table [x=h, y=eUCRMPS] {Plots/NavierStokes_num100ef0_2_7.txt};
	\addlegendentry{$\eps_0^{Mps}(\uvec_h)$}
	\logLogSlopeTriangle{0.9}{0.2}{0.23}{2}{black};
	\end{axis}       
	\end{tikzpicture}
 	\caption{$\eps_0(\uvec_h)$ for $(\uvec,p) \in $ \eqref{Test case Green-Taylor }. }
\label{fig:NSgreentaylorU}
\end{figure}

\begin{figure}[!ht]
	\centering
	\begin{tikzpicture}
	\begin{axis}[
	height = 7cm, 
	width = 10cm,  
	xlabel = {$h$},
	ylabel = {$\eps_0(p_h)$},
	ymajorgrids=true,
	major grid style={black!50},
	xmode=log, ymode=log,
	xmin=1e-2, xmax=1e-1,
	ymin=1e-2, ymax=1e1,
	yticklabel style={
		/pgf/number format/fixed,
		/pgf/number format/precision=0},
	scaled y ticks=false,
	xticklabel style={
		/pgf/number format/fixed,
		/pgf/number format/precision=0},
	scaled y ticks=false,  
	legend style={at={(1.3,0)},anchor=south east},
	legend columns=1
	],
	\addplot[color=mygreen,mark=triangle*] table [x=h, y=ePCRP0] {Plots/NavierStokes_num100ef0_2_7.txt};
	\addlegendentry{$\eps_0^{CR}(p_h) $ }
	\addplot[color=myred,mark=triangle*] table [x=h, y=ePCRP0P1] {Plots/NavierStokes_num100ef0_2_7.txt};
	\addlegendentry{$\eps_0^{\textrm{Trio}}(p_h)$}
	\addplot[color=myblue,mark=triangle*] table [x=h, y=ePCRMPS] {Plots/NavierStokes_num100ef0_2_7.txt};
	\addlegendentry{$\eps_0^{Mps}(p_h)$}
	\logLogSlopeTriangle{0.9}{0.2}{0.23}{1}{black};
	\end{axis}       
	\end{tikzpicture}
 	\caption{$\eps_0(p_h)$ for $(\uvec,p) \in $ \eqref{Test case Green-Taylor }.}
\label{fig:NSgreentaylorP}
\end{figure}
\newpage
 

\section{Conclusion and perspectives}

The purpose of this work is to present a new discretisation for the gradient of pressure. This scheme presents similar result to $\bP^1_{nc}-(P^0+P^1)$ discretization.   
However, some points have been left out of the scope of this work and deserve further investigation:
\begin{itemize}
	\item On the boundary, the continuity of the gradient flows can not be applied. We need boundary conditions to complete the system of elimination of the auxiliary unknowns \eqref{eq:continuitefacesnormales}. If the problem has, for the pressure:
	\begin{itemize}
		\item[-] Dirichlet boundary condition: we can evaluate the value of the auxiliary unknowns on the boundary.
		\item[-] Neumann boundary condition: we can evaluate the value of normal component of the pressure gradient on the boundary. 
	\end{itemize}
	Otherwise, we can keep the auxiliary unknowns and complete the problem with other equations.\\
	
	\item The section 3 shows that our scheme provides a benefit to the classic $\bP^1_{NC}-P^0$ discretisation but $\bP^1_{nc}-(P^0+P^1)$ has an additional superconvergence case. This property disappears in 3D, unless we add pressure degree of freedom on edges which turns out to be costly in computer memory. In that case, the MPFA scheme and $\bP^1_{nc}-(P^0+P^1)$ give comparable results but with a duality between scheme stencil and memory footprint. A study will be carried out to compare the efficiency of the two schemes.
	\item The scheme seems numerically stable but the inf-sup condition has still not been proven. 
	\item The scheme is currently in development in the CEA thermohydraulic code TrioCFD and its implementation will allow to realise more test.
	\item The FECC scheme is an other gradient discretization scheme, which has similar properties to the MPFA scheme and can handle more general meshes \cite{lepotier12}. The same approach can be used to develop a new scheme on polyhedral meshes for the Navier-Stokes problem.
	
\end{itemize}

\newpage
\printbibliography

\end{document}